\documentclass[11pt,reqno]{amsart}
\usepackage{tikz}
\usepackage{subfig}
\usepackage{epstopdf}
\usepackage{epsfig}
\usepackage{rotating}
\usepackage{hyperref}
\usepackage{cleveref}
\graphicspath{{../figures/},{"../../Code/hierarchical models/FiguresAdamNew/"},{"../../Code/hierarchical models/Type0/"}}
\usepackage{tikz}
\usetikzlibrary{shapes,arrows}
\tikzstyle{decision} = [diamond, draw, fill=blue!20, 
    text width=4.5em, text badly centered, node distance=3cm, inner sep=0pt]
\tikzstyle{block} = [rectangle, draw, fill=blue!20, 
    text width=5em, text centered, rounded corners, minimum height=4em]
\tikzstyle{line} = [draw, -latex']
\tikzstyle{cloud} = [draw, ellipse,fill=red!20, node distance=3cm,
    minimum height=2em]

\usetikzlibrary{positioning}
\tikzset{main node/.style={circle,fill=blue!20,draw,minimum size=1cm,inner sep=0pt},  }

\usepackage{soul}
\usepackage{color}
\usepackage{amsthm}
\newcommand{\JX}[1]{{\color{red}{$^{\text{JX}}$[#1]}}}

\newtheorem{conjecture}{Conjecture}
\usepackage{hyperref}
\usepackage{comment}
\usepackage{graphicx} \usepackage{placeins}
\usepackage{essay-def}
\usepackage[titletoc]{appendix}

\title[WIM on Gaussian mixture models]{Scaling Limits of the Wasserstein information matrix on Gaussian Mixture Models}
\author[Li]{Wuchen Li}
\email{wuchen@mailbox.sc.edu}
\address{Department of Mathematics, University of South Carolina, Columbia, SC 29208.}
\author[Zhao]{Jiaxi Zhao} 
\email{e0911835@u.nus.edu}
\address{Department of Mathematics, National University of Singapore.}
\thanks{W. Li's work is supported by AFOSR MURI FP 9550-18-1-502, AFOSR YIP award No. FA9550-23-1-0087, NSF RTG: 2038080, and NSF DMS-2245097.}
\begin{document}
\maketitle
\begin{abstract}
	We consider the Wasserstein metric on the Gaussian mixture models (GMMs), which is defined as the pullback of the full Wasserstein metric on the space of smooth probability distributions with finite second moment. It derives a class of Wasserstein metrics on probability simplices over one-dimensional bounded homogeneous lattices via a scaling limit of the Wasserstein metric on GMMs. Specifically, for a sequence of GMMs whose variances tend to zero, we prove that the limit of the Wasserstein metric exists after certain renormalization.
 Generalizations of this metric in general GMMs are established, including inhomogeneous lattice models whose lattice gaps are not the same, extended GMMs whose mean parameters of Gaussian components can also change, and the second-order metric containing high-order information of the scaling limit. We further study the Wasserstein gradient flows on GMMs for three typical functionals: potential, internal, and interaction energies. Numerical examples demonstrate the effectiveness of the proposed GMM models for approximating Wasserstein gradient flows.
\end{abstract}

\section{Introduction}

Optimal transport plays essential roles in mathematics, physics, and engineering \cite{Ambrosio2008Gradient, villani2021topics} with applications in data sciences \cite{peyre2019computational}. It studies a class of distance functionals, namely Wasserstein distances, between probability density distributions. Compared to classical divergences such as Kullback-Leibler divergence and Jensen–Shannon divergence, the Wasserstein distance has more robustness when dealing with distributions supported on low dimensional space. Thus, it is widely used in modern generative modeling and classification problems. Besides, gradient flows on Wasserstein spaces relate to problems in various fields, including partial differential equations \cite{JKO, Mielke_2011}, optimization problems in artificial intelligence \cite{NIPS2013_af21d0c9, gao2022master}, inverse problems \cite{Wangflow}, and geometric analysis \cite{LV, doi:10.1002/cpa.20060, e25050786}. 

Due to the increasing complexity in data science and machine learning \cite{IG, NG}, classical statistics and machine learning have provided several classes of nonlinear models in probability density spaces. One classical example is the Gaussian mixture model. Because of its strong power to approximate the multimodal distribution, the Gaussian mixture family has been widely applied to statistical estimation and sampling problems \cite{reynolds2009gaussian, IG}. One well-known fact is that the Fisher-Rao metric in a GMM, namely the ``Fisher information matrix'' \cite{IG, NG, ay2017information} (different from the Fisher information functional), has analytical properties. Unlike information geometry studies, understanding the Wasserstein metric in GMM remains an open question. Natural questions arise: {\em What are explicit formulations and scaling limits of the Wasserstein metrics restricted in Gaussian mixture models? Can we apply them to compute and approximate Wasserstein dynamics, such as Wasserstein gradient flows?}

In this paper, following the transport information geometry \cite{LiG, 10.1063/5.0012605, bregmanTIG}, we analyze the Wasserstein pullback metric in GMMs and define a Wasserstein metric on probability simplex as the scaling limit of metrics on a family of GMMs. We first introduce (homogeneous) Gaussian mixture models (GMMs), i.e.
\bequn
	\begin{aligned}
	\rho: \Theta & \rightarrow \mcP\lp \mbR \rp, \quad \Theta \subset \mbR^{N - 1},			\\
	\theta & \mapsto \rho_{\theta} = \sum_{i = 1}^{N - 1} \theta_i \lp \rho_{i + 1} - \rho_i \rp + \rho_1, \  1  = \theta_0 > \theta_1 > \cdots > \theta_{N - 1} > \theta_N = 0,			\\
	\rho_i & = \mcN\lp \mu_i, \sigma^2 \rp, \ \mu_1 < \mu_2 < \cdots < \mu_{N - 1} < \mu_N,
	\end{aligned}
\eequn
where $\mcP\lp \mbR \rp = \{ \rho \in C^{\infty}\lp \mbR \rp \big| \int_{\mbR} \rho\lp x \rp dx = 1\}$. In the preliminary stage, we consider the simplified model where the means and variances of Gaussian components are fixed to be constant while only mixing coefficients are allowed to vary. For more details in GMMs, c.f. \Cref{IM-GMM} and \cite{reynolds2009gaussian}. Denote the pull-back metric, i.e., Wasserstein information matrices (WIM) \cite{li2023wasserstein} on these models as $G_W\lp \theta; \sigma \rp = \lp G_W\lp \theta; \sigma \rp \rp_{1 \leq i, j \leq N-1} \in \mbR^{(N - 1) \times (N - 1)}$. The rigorous definition of WIM will be stated in \Cref{limit} and see more examples and properties in \cite{IG, NG, li2023wasserstein, LI2022109645, RePEc:spr:aistmt:v:74:y:2022:i:1:d:10.1007_s10463-021-00788-1}. As the variance parameter $\sigma$ tends to $0$, we have a weak convergence of Gaussian mixtures to Dirac mixtures. Hence, we define a scaling limit of WIMs under $\sigma \rightarrow 0$ to be a scaling Wasserstein metric on probability simplex.
\begin{Thm}[Informal Statement]\label{thm1-informal}
	For a 1D homogeneous GMM with differences between adjacent components given by $d$, a scaling limit of WIMs is given by
	\bequn
		\lim_{\sigma \rightarrow 0} \frac{\lp G_W\lp \theta; \sigma \rp \rp_{ij}}{K\lp \sigma \rp} = \frac{\delta_{ij}}{\sqrt{p_ip_{i + 1}}}, \quad K\lp \sigma \rp = \sqrt{2\pi^3}\frac{\sigma^3}{d}e^{\half \lp \frac{d}{2\sigma}\rp^2},
	\eequn
\end{Thm}
where $\delta_{ij}$s are the Kronecker symbols. \Cref{thm1-informal} is one of the main results of this paper. We also extend the results of scaling Wasserstein metrics to general 1D situations, including models with different gaps, higher-order expansions of metrics, as well as models with free mean parameters in \Cref{thm:inW-limit}, \Cref{thm:second-order-metric}, and \Cref{thm:extW-limit} respectively.

Specifically, we argue that \cref{equ:W-limit-matrix} is ``canonical'' compared to other definitions such as \cite{MAAS20112250, refId0, doi:10.1137/19M1243440}. We sketch the proof and calculation on Laplace mixture models in \Cref{L-limit} and leave the rigorous proof for general mixture models to future study. This establishes the universality of this WIM since it does not depend on specific models chosen to approximate the probability simplex. On the other hand, this definition is also important in defining and studying gradient flows on GMMs: since we obtain the scaling Wasserstein metric via dividing the metric tensor by a constant factor $K\lp \sigma \rp$, it is not difficult to find out that gradient flows relate to ordinary one by a time reparameterization, i.e., time scaling.

Studying the optimal transport theory in GMMs has attracted researchers' attention for a long time. In \cite{Malago2}, the author studies the properties of Wasserstein gradient flows in the mixture geometry of a discrete exponential family. The work in \cite{chen2018optimal} is mostly close to us; they defined the Wasserstein distance over GMMs via the original coupling definition of the transportation cost, which is parallel to our definition via pullback metric. They studied metric properties and geodesics of this metric and compared it with the original Wasserstein-2 distance. Moreover, numerical experiments are conducted on barycenters in GMMs to visualize the geodesics.

As a graphical model, GMM can also be analyzed through the lens of probability simplex \cite{chen2018optimal}. In literature, researchers have attempted to establish an optimal transport theory on graphs \cite{MAAS20112250}. \cite{doi:10.1137/19M1243440} studies the scaling limit of the Wasserstein metric on graphs using the Benamou-Brenier formulation and establishes the lower and upper bounds of the limiting metric. In \cite{MAAS20112250}, the author defines a metric over the probability simplex associated with a fixed Markov chain. Under this metric, the Markov chain can be identified with he gradient flow of entropy. Furthermore, in \cite{GLADBACH2020204}, the same authors prove that in a 1D periodic setting, discrete transport metrics
converge to a limiting transport metric with a non-trivial effective mobility. Compared to previous work, we focus on the 1D models and use the closed-form formula for the Wasserstein metric in 1D \cite{villani2021topics} to define the WIMs for GMMs. Therefore, we can obtain more quantitative information of the Wasserstein metric on these models and at the same time prove the existence of the scaling limit of Wasserstein metric.

The paper is organized as follows. In \Cref{limit}, we establish the scaling limits of Fisher and Wasserstein information matrices. Next, we generalize results to general models, such as inhomogeneous GMMs, higher-order GMM metrics, and extended GMMs in \Cref{generalize}. Last, in \Cref{gradient-flow}, we derive the gradient flows on these scaling Wasserstein metrics and relate them to gradient flows on density manifold. Numerical experiments and further discussions are left to \Cref{experiment} and \Cref{discussion}, respectively.

\section{A scaling limit of Fisher-Rao and Wasserstein metrics in Gaussian mixture models}\label{limit}

In this section, we study a scaling limit of Fisher and Wasserstein metrics on 1D GMMs. They are defined via the pullback of the metric on the density manifold. The scaling limit is defined as the limit of variances of Gaussian components tend to $0$.
\subsection{Information matrices and Gaussian mixture models}\label{IM-GMM}

Let $\mathcal{X}\subset \mathbb{R}^n$ be a compact manifold with a smooth Riemannian structure. This can be the computational domain of a scientific application or the sample space associated with statistical inference problems. Let $\mathcal{P}^{2, \infty}(\mathcal{X})$ denote the space of probability distributions over $\mathcal{X}$ which are absolutely continuous w.r.t. the Lebesgue measure on $\mathcal{X}$ and have smooth positive density functions. Moreover, one further requires these distributions to have bounded second moments, which can be calculated via the embedding of $\mathcal{X}$ into $\mathbb{R}^n$. For those who do not have much familiarity with the Riemannian geometry, $\mathcal{X}$ can be simply regarded as the whole Euclidean space $\mathbb{R}^n$. Given a metric tensor $g$ on $\mathcal{P}^{2, \infty}(\mathcal{X})$, we call $(\mathcal{P}^{2, \infty}(\mathcal{X}),g)$ density manifold \cite{Lottgeo, kriegl1997convenient}. Consider a parameter space $\Theta\subset \mathbb{R}^d$ and a parameterization function
\begin{equation*}
    \rho \colon \Theta \rightarrow \mathcal{P}^{2, \infty}( \mcX ), \quad \theta \mapsto \rho_{\theta}(x)
\end{equation*}
which can also be viewed as
$\rho \colon \mathcal{X}\times \Theta\rightarrow \mathbb{R}$. We assume the distributions have smooth density functions. The image of $\Theta$ under the mapping $\rho$, i.e., $\rho(\Theta)$ is named a statistical model. Suppose that $f, g$ are two vector-valued functions on $\mathcal{X}$, we denote $\langle f, h\rangle=\int_{\mathcal{X}} (f(x), h(x)) dx$ as the $L^2(\mathcal{X})$ inner product, where $dx$ refers to the Lebesgue measure on $\mathcal{X}$. We denote $\lp v,w \rpz= v \cdot w$ as the (pointwise) Euclidean inner product of two vectors. On an arbitrary parametric space $\rho: \Theta \rightarrow \mcP\lp \mcX \rp$, we define a pull-back metric $g^{\Theta}$ and an information matrix $G$ associated to a metric $g$ on density manifold $\mcP\lp \mcX \rp$. See a related study in \cite{li2023wasserstein}.
\begin{Def}[Pull-back metric \& information matrix]
	Consider density manifold $(\mathcal{P}(\mcX), g)$ with a metric tensor $g$, i.e. $\forall \nu \in \mcP\lp \mcX \rp, g\lp \nu \rp: T_{\nu}\mcP\lp \mcX \rp \times T_{\nu}\mcP\lp \mcX \rp \rightarrow \mbR$ bilinear, and a smoothly parametrized space $\rho: \Theta \rightarrow \mcP\lp \mcX \rp$ with parameter $\theta\in\Theta \subset \mbR^d$. 
	Then the pull-back metric $g^{\Theta}$ of $g$ onto this parameter space $\Theta$ is given by
	\begin{equation*}
		\begin{aligned}
			g^{\Theta}\lp \rho_{\theta} \rp & : T_{\rho_{\theta}}\Theta \times T_{\rho_{\theta}}\Theta \rightarrow \mbR		\\
			g^{\Theta}\lp \rho_{\theta} \rp\lp v_1, v_2 \rp & = g\lp \rho_{\theta} \rp\lp v_1, v_2 \rp, \forall v_1, v_2 \in T_{\rho_{\theta}}\Theta \subset T_{\rho_{\theta}}\mcP\lp \mcX \rp.
		\end{aligned}
	\end{equation*}
	Denote the information matrix associated with this model as $G\lp \theta \rp \subset \mbR^{d \times d}, \forall \theta \in \Theta$
	\begin{equation*}
		G(\theta)_{ij} = g^{\Theta}\lp \rho_{\theta} \rp\lp \p_{\theta_i}\rho_{\theta}, \p_{\theta_j}\rho_{\theta} \rp.
	\end{equation*}
\end{Def}

In this paper, we restrict our attention to the special one-dimensional sample space $\mcX = \mbR$. Consider an 1D GMM $\rho: \Theta \rightarrow \mcP\lp \mbR \rp, \Theta \subset \mbR^{N - 1}$ with
\bequ\label{equ:GMM}
	\theta \mapsto \rho_{\theta} = \sum_{i = 1}^{N - 1} \theta_i \lp \rho_{i + 1} - \rho_i \rp + \rho_1, \qquad 1 > \theta_1 > \cdots > \theta_{N-1} > 0, \quad \forall i = 1,...,N - 1,\tag{GMM}
\eequ
where each component is fixed to be a Gaussian $\rho_i(x) = \frac{1}{\sqrt{2\pi}\sigma_i}e^{-\frac{(x-\mu_i)^2}{2\sigma_i^2}}, i = 1,2,...,N$ and components are listed in their means' increasing order $\mu_1 < \mu_2 < ... < \mu_{N-1} < \mu_{N}$. As mentioned before, we consider the simplified model where only mixing coefficients are allowed to vary. In \Cref{sec:ext-GMM}, we will relax this constraint and consider GMM with more degrees of freedom. We further postulate $\theta_0 = 1, \theta_N = 0$. The parameters $\theta_i$s do not have any probabilistic meaning, we thus introduce another group of coordinate variables $p_i, i = 1,2...N$ for GMMs, i.e.
\bequn
	\sum_{i = 1}^{N - 1} \theta_i \lp \rho_{i + 1} - \rho_i \rp + \rho_1 = \sum_{i = 1}^{N} \lp \theta_{i - 1} - \theta_i \rp \rho_i = \sum_{i = 1}^{N} p_i \rho_i.
\eequn
Writing in this form, the GMM has close relation with the probability simplex, i.e. $\lbb p_i, i\in [N] \rbb$ represents a point in the probability simplex. \cite{chen2018optimal} use heavily this connection to study the properties of optimal transport on GMMs. Consequently, any metric defined on GMMs can also be viewed as a metric on the probability simplex. Throughout the paper, we will use the coordinates $\theta_i$s and $p_i$s interchangeably as $\theta$-coordinate simplifies the presentation of the WIM while $p$-coordinate is better for understanding. To simplify the deduction, we postulate a homogeneous assumption, i.e. \Cref{assm} of the model in this section and name it 1D homogeneous GMM.
\begin{Asm}[Homogeneity]\label{assm}
	The standard variances of all the components coincide and the differences between adjacent means of Gaussian components are the same:
	\bequn
		\begin{aligned}
			\sigma_i = \sigma,\qquad i = 1,2,...,N-1,N,		\\
			\mu_i - \mu_{i-1} = d,\qquad i = 2,3,...,N-1,N.
		\end{aligned}
	\eequn
\end{Asm} 
In following examples, we investigate the Fisher and Wasserstein information matrices in the Gaussian mixture models.
\begin{Ex}[Fisher geometry of Gaussian mixture model]
One can identify the tangent space of density manifold $\mcP$ at arbitrary distribution $\rho$ with $C_0^{\infty}$:
\begin{equation}\label{tangent-space} 
	\begin{aligned}
    		T_{\rho}\mcP(\mcX) & \simeq C_0^{\infty}(\mathcal{X}) = \{f \in C^{\infty}(\mcX) | \int_{\mcX}f(x)dx = 0\},
	\end{aligned}
\end{equation}
where $C^{\infty}(\mcX)$ is the space of continuous functions on $\mcX$. Then, the Fisher-Rao metric is given as
\bequn
	\begin{aligned}
		(\text{Fisher-Rao}): & \ 	g_F\lp v_1, v_2 \rp = \int_{\mcX} \frac{v_1(x) v_2(x)}{\rho(x)}dx, \quad v_1, v_2 \in T_{\rho}\mcP\lp \mcX \rp,
	\end{aligned}
\eequn
where we omit the dependence of $g_F\lp \rho \rp$ on $\rho$. For a 1D GMM \cref{para}, Fisher information matrices are provided as
\bequ\label{equ:FIM-GMM}
	\begin{aligned}
		\lp G_F \lp \theta \rp\rp_{ij} &  = \int_{\mbR} \frac{\p_{\theta_i}\rho_{\theta}(x) \p_{\theta_j} \rho_{\theta}(x) }{\rho_{\theta}(x)} dx= \int_{\mbR} \frac{\lp \rho_{i + 1}(x) - \rho_i(x) \rp \lp \rho_{j + 1}(x) - \rho_j(x) \rp}{\rho_{\theta}(x)} dx.
	\end{aligned}
\eequ
\end{Ex}
\begin{Ex}[Wasserstein geometry of Gaussian mixture model]
Based on the identification \cref{tangent-space} of the tangent space in last example, we have
\bequn
	\begin{aligned}
		(\text{Wasserstein}): & \ 	g_W\lp v_1, v_2 \rp = \int_{\mcX} v_1(x) \Phi_2(x) dx, \ v_2(x) + \nabla \cdot \lp \rho(x) \nabla \Phi_2(x) \rp = 0, \quad v_1, v_2 \in T_{\rho}\mcP\lp \mcX \rp,
	\end{aligned}
\eequn
Moreover, when we focus on 1D cases, i.e. $\mcX = \mbR$, the metric for Wasserstein geometry has a closed-form solution, namely
\bequn
	g_W\lp v_1, v_2 \rp = \int_{\mbR} \frac{F_1(x) F_2(x)}{\rho(x)} dx, \quad F_i\lp x \rp = \int_{-\infty}^{x} v_i\lp s \rp ds, i = 1, 2.
\eequn
Hence, for a 1D GMM \cref{para}, the WIMs are provided as
\bequ\label{equ:WIM-GMM}
	\begin{aligned}
		\lp G_W\lp \theta \rp \rp_{ij} & = \int_{\mbR} \frac{\p_{\theta_i}F_{\theta}(x) \p_{\theta_j} F_{\theta}(x) }{\rho_{\theta}(x)} dx= \int_{\mbR} \frac{\lp F_{i + 1}(x) - F_i(x) \rp \lp F_{j + 1}(x) - F_j(x) \rp}{\rho_{\theta}(x)} dx,
	\end{aligned}
\eequ
where $F_i$ is the cumulative distribution function of density $\rho_i$.

\end{Ex}
\subsection{The scaling limit of information matrices}
Although the analytic formula for the WIM in 1D GMM exists, the integral \cref{equ:WIM-GMM} is impossible to evaluate explicitly. Thus, we consider the following scaling limit.
\begin{Def}[Scaling limit of the information matrices in GMMs]\label{def:scaling}
    Let $\Theta(\sigma^2) \subset \mcP(\mbR)$ denote a family of GMMs indexed by the variance of Gaussian components, i.e., they share the same means for each component but different variance $\sigma^2$:
\bequ\label{para}
	\theta \mapsto \rho(\cdot|\theta, \sigma^2) = \sum_{i = 1}^{N - 1} \theta_i \lp \rho(\cdot|\mu_{i+1}, \sigma^2) - \rho(\cdot|\mu_i, \sigma^2) \rp + \rho(\cdot|\mu_1, \sigma^2), 
\eequ
where $\rho(\cdot|\mu_i, \sigma^2) = \mcN(\mu_i, \sigma^2)$. Denote $G\lp \theta; \sigma \rp$ as the information matrices associated to the model $\Theta(\sigma^2)$. In our context, this can be either Fisher $G_F\lp \theta; \sigma \rp$ or Wasserstein $G_W\lp \theta; \sigma \rp$. Consider the limit as $\sigma^2 \rightarrow 0$, if there exists a function $K(\sigma)$ of $\sigma$ such that the following limit exists
\bequ\label{equ:scaling-limit}
    \lim_{\sigma \rightarrow 0} \frac{G\lp \theta; \sigma \rp}{K(\sigma)} = \wtd G\lp \theta \rp.
\eequ    
We call the limit matrix $\wtd G\lp \theta \rp$ the scaling limit of the information matrices in GMMs, or briefly scaling information matrices.
\end{Def}
There are two reasons to consider this scaling limit. The first one is that both Fisher and Wasserstein geometry possess a well-defined scaling limit, which is one of the main discoveries of this work. Another motivation is the following: as the standard variance $\sigma$ of Gaussian $\mcN\lp \mu, \sigma^2 \rp$ tends to zero, the distribution converges weakly to a Dirac measure $\delta_{\mu}$ centered at its mean $\mu$. Specifically, Gaussian mixture models converge weakly to Dirac mixture models. Consequently, this limiting behavior of the information matrices characterizes the corresponding geometry on the Dirac mixture model, which can also be understood as probability simplices or discrete graphical models. If the limit exists, it is a candidate of the Wasserstein metric on probability simplex, which is important in both communities of optimal transport and graph theory \cite{MAAS20112250, li2022mean, slepcev2022nonlocal, gao2022master}. In comparison, our scaling metrics \Cref{def:scaling} differs from \cite{MAAS20112250} in the sense that metrics $\wtd G\lp \theta \rp$ depends only on the structure of graphs, while theirs depends also on transition kernels of jump processes. Moreover, scaling limit of the Wasserstein geometry is also studied in \cite{GLADBACH2020204}, but it considers the opposite direction. While they explore the convergence of optimal transport on grid to continuous situation as the size of grid tends to $0$, we define the optimal transport on graphs by approximating it via continuous models.

As a warm-up, we calculate the scaling limit of the Fisher-Rao metric, which is simple in the sense that one can directly set the scaling factor $K(\sigma)$ to $1$ to obtain the desired limit.
\begin{Thm}
	For a 1D homogeneous GMM, the scaling limit of Fisher information matrices is given by
	\bequ\label{Fisher-limit}
		\lim_{\sigma \rightarrow 0} \lp G_F\lp \theta; \sigma \rp \rp_{ij} = \lbb 
			\begin{aligned}
				& \frac{1}{p_i} + \frac{1}{p_{i + 1}}, \quad i = j,	\\
				& - \frac{1}{p_i}, \qquad \quad i = j + 1,			\\
				& - \frac{1}{p_{i + 1}}, \quad  \quad  j = i + 1,		\\
				& 0, \qquad \qquad  \ \text{ otherwise,}
			\end{aligned}\rpt\tag{F-limit}
	\eequ
	or in matrix form
	\bequ\label{equ:Fisher-limit-matrix}
		G_{\wtd F}\lp \theta \rp = \lim_{\sigma \rightarrow 0} G_F\lp \theta; \sigma \rp = \begin{pmatrix}
			\frac{1}{p_1} + \frac{1}{p_{2}} & - \frac{1}{p_2} & 0 & \cdots & 0 & 0 			\\
			- \frac{1}{p_2} & \frac{1}{p_2} + \frac{1}{p_{3}} & - \frac{1}{p_3} & \cdots & 0 & 0	\\
			0 & - \frac{1}{p_3} & \frac{1}{p_3} + \frac{1}{p_{4}} & \cdots & 0 & 0			\\
			\vdots & \vdots & \vdots & \ddots & \vdots & \vdots 						\\
			0 & 0 & 0 & \cdots & -\frac{1}{p_{N - 1}} & \frac{1}{p_{N - 1}} + \frac{1}{p_{N}}
		\end{pmatrix}.\tag{F-limit-matrix}
	\eequ
\end{Thm}
\begin{proof}
	Recall that for a 1D homogeneous GMM, the Fisher information matrix reduces to \cref{equ:FIM-GMM}. For simplification, we omit the dependence of $\rho_{\theta}\lp x \rp$ on $\sigma$. Readers should not be confused when we take the limit of the above quantity as $\sigma$ tends to $0$. Thus, it suffices to prove the following relation
	\bequn
		\lim_{\sigma \rightarrow 0}\int \frac{\rho_i\lp x \rp \rho_j\lp x \rp}{\rho_{\theta}\lp x \rp}dx = \frac{\delta_{ij}}{p_i},
	\eequn
	where $\delta_{ij}$s are Kronecker symbols. Since we know that a sequence of Gaussians with same mean $\mu_i$ and variances shrink to $0$ converges weakly to the Dirac measure at $\mu_i$, i.e. $\mcN\lp \mu_i, \sigma \rp \stackrel{w}{\longrightarrow} \delta_{\mu_i}$. We have
	\bequn
		\begin{aligned}
			& \ \lim_{\sigma \rightarrow 0}\int \frac{\rho_i\lp x \rp \rho_j\lp x \rp}{\rho_{\theta}\lp x \rp}dx		\\
			= & \ \lim_{\sigma \rightarrow 0}\mbE_{X \sim \rho_i} \frac{\rho_j\lp X \rp}{\rho_{\theta}\lp X \rp}	\\
			= & \ \mbE_{X \sim \delta_{\mu_i}} \lim_{\sigma \rightarrow 0}\frac{\rho_j\lp X \rp}{\rho_{\theta}\lp X \rp}		\\
			= & \lim_{\sigma \rightarrow 0}\frac{\rho_j\lp \mu_i \rp}{\rho_{\theta}\lp \mu_i \rp} = \frac{\delta_{ij}}{p_i}.
		\end{aligned}
	\eequn
 We exchange the expectation and limit as the integrand $\frac{\rho_j\lp X \rp}{\rho_{\theta}\lp X \rp}$ is bounded from above by $\frac{1}{p_i}$ and thus uniformly integrable.
\end{proof}
\begin{Rem}
	Readers who are familiar with Fisher information geometry can realize that the scaling limit \cref{equ:Fisher-limit-matrix} obtained is exactly the Fisher information matrix on the probability simplex 
	\bequn
		\circ-\circ-\circ- \cdots -\circ-\circ,
	\eequn
	under $\theta$-parameterization \cref{para}, except that Gaussian components are replaced by Dirac measures on each node. In other word, the scaling limit of Fisher-Rao metric on continuous models coincides with Fisher geometry on discrete models. This identification in Fisher geometry implies that the definition via scaling limit is at least canonical in Fisher geometry. This motivates us to study the counterpart in Wasserstein geometry in the next subsection.
\end{Rem}
\subsection{The scaling limit of Wasserstein metric}

In this section, we study the scaling limits of WIMs $G_{\wtd W}$. We put this derivation into general scope by considering not only the scaling behavior over GMMs, but also all the mixture models of the form $\rho\lp x; \theta \rp = \sum_{i = 1}^N p_i \rho\lp x - \mu_i \rp$ under scaling $\rho \lp x;\sigma \rp = \sigma\rho \lp \sigma x \rp$, where $\rho(\cdot)$ is a probability density. Specifically, we will perform detailed calculation over Guassian and Laplace families and remark on general families.

\begin{figure}[ht]\label{fig:pdf-cdf}
  \centering
  \label{mixture}\centerline{\includegraphics[width=0.52\linewidth]{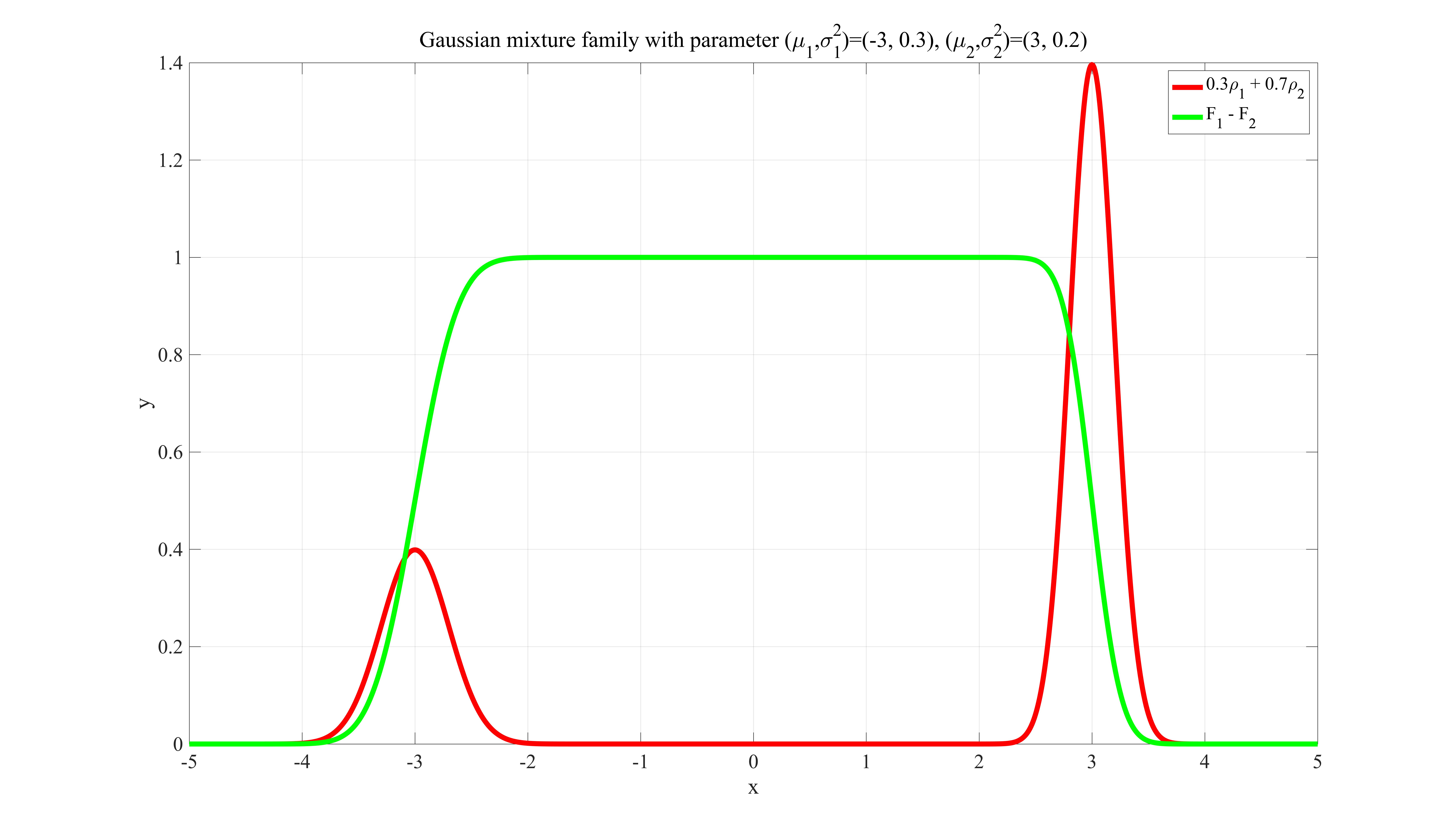}}
  \caption{This figure plots an example of the function $\p_{\theta_i} F_{\theta}\lp x \rp$ for a GMM.}
\end{figure}

Recall that for 1D models the Wasserstein inner product is given by \cref{equ:WIM-GMM}. However, in mixture models, the term $\p_{\theta_i} F_{\theta}\lp x \rp = F_{i + 1}\lp x \rp - F_{i}\lp x \rp$ has the following behaviors shown in \Cref{fig:pdf-cdf}.
We plot the density function of a 2-mixture model and $F_{2}\lp x \rp - F_{1}\lp x \rp$. It is easy to observe that on a large part of the interval $\lb \mu_1, \mu_2 \rb$, the function $\p_{\theta_i} F_{\theta}\lp x \rp$ stays close to $1$ while $\rho_{\theta}\lp x \rp$ is relatively small, as $\sigma$ tends to $0$. Consequently, as the scaling parameter $\sigma$ tends to $0$, \cref{equ:WIM-GMM} blows up, i.e. $\lim_{\sigma \rightarrow 0}G_W\lp \theta; \sigma \rp$ does not exist. This indicates the importance of considering the scaling limit \cref{equ:scaling-limit} in \Cref{def:scaling}. And the primary goal becomes quantifying the scaling factor $K(\sigma)$. The main result is as follows
\begin{Thm}\label{thm:W-limit}
	For a 1D homogeneous GMM with difference between adjacent components given by $d$, the scaling limit of WIMs is given by
	\bequ\label{equ:W-limit}
		\lim_{\sigma \rightarrow 0} \frac{\lp G_W\lp \theta; \sigma \rp \rp_{ij}}{K\lp \sigma \rp} = \frac{\delta_{ij}}{\sqrt{p_ip_{i + 1}}}, \quad K\lp \sigma \rp = \sqrt{2\pi^3}\frac{\sigma^3}{d}e^{\frac{d^2}{8\sigma^2}}.\tag{naive-W-limit}
	\eequ
	In matrix form,
	\bequ\label{equ:W-limit-matrix}
		G_{\wtd W}\lp \theta \rp = \lim_{\sigma \rightarrow 0} \frac{G_W\lp \theta; \sigma \rp}{K\lp \sigma \rp} = \begin{pmatrix}
			\frac{1}{\sqrt{p_1p_{2}}} & 0 & 0 & \cdots & 0 & 0 			\\
			0 & \frac{1}{\sqrt{p_2p_{3}}} & 0 & \cdots & 0 & 0	\\
			0 & 0 & \frac{1}{\sqrt{p_3p_{4}}} & \cdots & 0 & 0			\\
			\vdots & \vdots & \vdots & \ddots & \vdots & \vdots 						\\
			0 & 0 & 0 & \cdots & 0 & \frac{1}{\sqrt{p_{N - 1}p_{N}}}
		\end{pmatrix}\tag{naive-W-limit-matrix}.
	\eequ
\end{Thm}
\begin{Rem}
	Comparing \cref{equ:W-limit-matrix} and \cref{equ:Fisher-limit-matrix}, the key difference is that scaling limit of WIM requires an additional rescaling factor $K(\sigma)$. The reason is that Fisher geometry is invariant under parameterization (indeed, Fisher-Rao metric can be characterized as the unique metric satisfies a kind of invariance under the variable transformation), but Wasserstein geometry clearly does not satisfy the same property.
\end{Rem}

The rest of this section is devoted to the proof of \Cref{thm:W-limit}, which contains a detailed analysis using Laplace's method. Before diving into technical proofs, We provide a brief overview on the proof idea. We analyze the WIM via the following reduction
\bequ\label{equ:reduction}
	\begin{aligned}
		\lp G_W\lp \theta \rp \rp_{ii} = \int_{\mbR} \frac{\lp F_i - F_{i + 1} \rp^2}{\rho_{\theta}}dx \stackrel{\Delta_1}{\Longrightarrow} \int_{\mu_{i}}^{\mu_{i+1}} \frac{dx}{p_{i}\rho_{i} + p_{i+1} \rho_{i+1} } \stackrel{\Delta_2}{\Longrightarrow}  \text{Laplace's method}.
	\end{aligned}
\eequ
In the $\Delta_1$ reduction, we restrict the integration domain from $\mbR$ to $[\mu_{i}, \mu_{i+1}]$ which contains the most significant part of the numerator $F_i(x) - F_{i+1}(x)$, c.f. \Cref{fig:pdf-cdf}. Meanwhile, we simplify the total density $\rho_{\theta}$ to $p_{i}\rho_{i} + p_{i+1} \rho_{i+1}$ which are the two nearest components and $(F_i(x) - F_{i+1}(x))^2$ to $1$. In the next $\Delta_2$ reduction step, we perform the limit $\sigma \rightarrow 0$ and deduce the asymptotic formula for the integral using a variant of the Laplace's method \cite{murray2012asymptotic}. 

The rigorous derivation of $\Delta_1, \Delta_2$ reduction are summarized in the following two propositions.
\begin{Prop}[$\Delta_1$ reduction]\label{prop:Delta-1}
	Consider the WIM on the homogeneous GMM \Cref{assm} with variance $\sigma^2$ and gap $d$. The non-diagonal term has the following upper bound:
	\begin{equation}
		\begin{aligned}
			0 \leq \lp G_W \rp_{ij} & \leq \frac{3d\sigma^2 M}{\min_i p_i}\lp 1 + \frac{2\sigma^2}{3d^2}e^{-d^2/2\sigma^2}\rp, \quad i \neq j,     \\
		\end{aligned}
	\end{equation}
 and the diagonal term can be bounded as
 \begin{equation}
		\begin{aligned}
			 \frac{\sigma^4M^2}{2\max\{p_i, p_{i+1}\}}+\lp 1 - \sqrt{\frac{2\sigma}{\pi}}e^{-1/2\sigma} - \frac{e^{-d^2/2\sigma^2}}{\min\{p_i, p_{i+1}\}} \rp I_i \leq \lp G_W \rp_{ii} \leq \frac{\sigma^4M^2}{2\min\{p_i, p_{i+1}\}} + I_i,     \\
		\end{aligned}
\end{equation}
with the integral $I_i$ given by
	 \begin{equation}\label{equ:Delta-2-integral}
	     I_i = \int_{\mu_{i}}^{\mu_{i + 1}} \frac{dx}{p_{i + 1}\rho_{i + 1} + p_{i} \rho_i }, \quad i \in [N-1].
	 \end{equation}
\end{Prop}
A glimpse over the bounds of non-diagonal terms and diagonal terms indicate two classed of matrix elements have significant scale difference. While the non-diagonal terms tend to $0$ with speed $O(\sigma^2)$ as $\sigma \rightarrow 0$, the diagonal terms tend to $\infty$ similar to the asymptotic integral $I_i$, which will be proved to with the order $e^{d^2/8\sigma^2}$. In a word, the WIM rapidly becomes diagonally dominant as $\sigma \rightarrow 0$, which is aligned with the conclusion in \Cref{thm:W-limit}.

To deal with the reduction $\Delta_2$, we have the following variant of Laplace's method. To maintain simplification, we focus on two specific density families, namely Gaussian and Laplace, to obtain their asymptotic analysis. Meanwhile, we establish a slightly more general version which can also handle a generalized model in \Cref{generalize}. 
\begin{Prop}[$\Delta_2$ reduction]\label{prop:Delta-2}
	We have the following asymptotics of the \cref{equ:Delta-2-integral}
	\bequ
		\begin{aligned}
			\lim_{\sigma \rightarrow 0} \frac{\int_0^{d} \frac{dx}{p_{i}\rho \lp x; \sigma \rp + p_{i+1} \rho 				\lp d - x; k\sigma \rp} }{\frac{\sqrt{2\pi} \sigma^3 }{p_{i}l}e^{\half \lp l/\sigma \rp^2}} & = g \lp k \rp,  \quad (\text{Gaussian})       \\
   \lim_{\sigma \rightarrow 0} \frac{\int_0^{d} \frac{dx}{p_{i}\rho \lp x; \sigma \rp + p_{i+1} \rho 				\lp d - x; k\sigma \rp} }{\frac{2\sigma^2}{p_i}e^{l/\sigma}} & = g \lp k \rp,  \quad (\text{Laplace})
		\end{aligned}
	\eequ
 where $g \lp k \rp = \int_0^{\infty} \frac{dy}{ 1 + y^{(k + 1)/k}	}$ and the parameter $l$ satisfies the following matching condition
	\bequ\label{equ:match-cond-prop}
		p_{i}\rho \lp l; \sigma \rp = p_{i+1} \rho \lp 	d - l; k\sigma \rp.
	\eequ
\end{Prop}

\begin{proof}[Proof of \Cref{thm:W-limit}]
It suffices to prove the result \cref{equ:W-limit} elementwise. For crossing terms $\lp G_W \rp_{ij}, i \neq j$, by \Cref{prop:Delta-1}, we conclude
\bequn
	\lim_{\sigma \rightarrow 0} \frac{\lp G_W \rp_{ij}}{K\lp \sigma \rp} \leq \lim_{\sigma \rightarrow 0} \frac{\frac{3d\sigma^2 M}{\min_i p_i}\lp 1 + \frac{2\sigma^2}{3d^2}e^{-d^2/2\sigma^2}\rp}{K\lp \sigma \rp} = 0.
\eequn
For diagonal terms, similar reduction via \Cref{prop:Delta-1} deduces
\bequn
	\lim_{\sigma \rightarrow 0} \frac{\lp G_W \rp_{ii}}{K\lp \sigma \rp} = \lim_{\sigma \rightarrow 0} \frac{\int_{\mu_{i}}^{\mu_{i+1}} \frac{dx}{p_{i}\rho_{i} + p_{i + 1} \rho_{i + 1} }}{K\lp \sigma \rp}.
\eequn
Applying \Cref{prop:Delta-2} to the case $k=1$, we obtain
\bequ
	\lim_{\sigma \rightarrow 0} \frac{\int_{\mu_i}^{\mu_{i+1}} \frac{dx}{p_{i}\rho_i + p_{i+1} \rho_{i+1}}}{\frac{\sqrt{2\pi} \sigma^3 }{p_{i}l}e^{\half \lp l/\sigma \rp^2}} = g \lp 1 \rp = \frac{\pi}{2},
\eequ
where $l$ is given by \cref{equ:l-gaussian}. Using this, we obtain
\begin{equation}
    \frac{\sqrt{2\pi} \sigma^3 }{p_{i}l}e^{\half \lp l/\sigma \rp^2} = \frac{2\sqrt{2\pi} \sigma^3 }{p_{i}d}\sqrt{\frac{p_i}{p_{i+1}}}e^{\half \lp d/2\sigma \rp^2} + O(\sigma^2).
\end{equation}
Combining with previous results we conclude the proof.
\end{proof}
\begin{Rem}[Scaling limit of WIMs in Laplace mixture models]\label{L-limit}
For Laplace family, the general workflow remains the same. Although the estimation for reduction $\Delta_1$ is not proved for this case (\Cref{prop:Delta-1} is only for GMMs), we comment here that correspondent bounds can be obtained in this family without much efforts. Applying \Cref{prop:Delta-2} to the case $k=1$, we obtain
\bequ
	\lim_{\sigma \rightarrow 0} \frac{\int_{\mu_i}^{\mu_{i+1}} \frac{dx}{p_{i}\rho_i + p_{i+1} \rho_{i+1}}}{\frac{2\sigma^2}{p_{i}}e^{l/\sigma}} = g \lp 1 \rp = \frac{\pi}{2},
\eequ
where $l$ is given as follows using the matching condition \cref{equ:match-cond}.
\begin{equation}\label{equ:l-laplace}
    \begin{aligned}
            l & = \frac{d}{k + 1} + \frac{k\sigma}{k+1}\log \frac{kp_{i}}{p_{i+1}}.
    \end{aligned}
\end{equation}
Via elementary calculation, we derive
\bequn
	\frac{2\sigma^2}{p_{i}}e^{l/\sigma} = \frac{2\sigma^2}{p_i} \lp \frac{p_ik}{p_{i+1}}\rp^{k/(k+1)}e^{\frac{d}{(k+1)\sigma}} = \frac{2\sigma^2}{\sqrt{p_ip_{i+1}}} e^{\frac{d}{2\sigma}}.
\eequn
By choosing $K\lp \sigma \rp = \pi\sigma^2 e^{\frac{d}{2\sigma}}$, we conclude that the scaling limit of $G_W\lp \theta; \sigma \rp$ in Laplace mixture models is given by
\bequn
	\lim_{\sigma \rightarrow 0} \frac{\lp G_W\lp \theta; \sigma \rp \rp_{ii}}{ \pi \sigma^2 e^{\frac{d}{2\sigma}}} = \frac{1}{\sqrt{p_ip_{i + 1}}},
\eequn
which coincides with $\cref{equ:W-limit-matrix}$ obtained from GMMs.
\end{Rem}
\begin{proof}[Proof of \Cref{prop:Delta-1}]
Throughout the proof we will use the following basic relation between the probability density function and cumulative density function $\Phi(x)$ of standard Gaussian distribution:
\begin{equation}
    \Phi(-x), 1-\Phi(x) \leq \frac{1}{ \sqrt{2\pi}x} e^{-x^2/2}, \quad \forall x>0.
\end{equation}
Specifically, we obtain a constant without dependence on $x$ by dividing into two cases $\lv x \rv \leq 1, \lv x \rv > 1$. We have
\begin{equation}
    \begin{aligned}
        \frac{\max\{\Phi(-x), 1-\Phi(x)\}}{\rho(x)} & < 1, \quad x > 1,        \\
        \frac{\max\{\Phi(-x), 1-\Phi(x)\}}{\rho(x)} & \leq \frac{\half}{\frac{1}{ \sqrt{2\pi}} e^{-1^2/2}} = \frac{\sqrt{2\pi e}}{2} := M, \quad 0 < x < 1.
    \end{aligned}
\end{equation}
This concludes that the constant $M$ can bound the LHS over $\mbR$. For general Gaussian with zero mean and $\sigma^2$ variance, we have for $x < 0$
\begin{equation}
    \Phi(x; \sigma) = \Phi(x/\sigma) \leq \frac{1}{\sqrt{2\pi}}\frac{\sigma}{\lv x \rv}e^{-x^2/2\sigma^2},
\end{equation}
which derives the following bound over $\mbR$:
\begin{equation}
    \begin{aligned}
        \frac{\max\{\Phi(-x;\sigma), 1-\Phi(x;\sigma)\}}{\rho(x;\sigma)} & < \sigma^2 M, \quad x \geq 0.        \\
    \end{aligned}
\end{equation}
Now, we can use bound the matrix element of WIM using this inequality.

\textit{Bound the off-diagonal terms}:
The integrand of the off-diagonal term has the following form:
\begin{equation}
    g_{ij}(x) := \frac{\lp F_{i}(x) - F_{i+1}(x) \rp \lp F_{j}(x) - F_{j+1}(x) \rp}{\rho_{\theta}(x)}, \quad i < j.
\end{equation}
Using the fact that
\begin{equation}
    \begin{aligned}
        & \lv F_{i}(x) - F_{i+1}(x) \rv < 1 \quad x \in \mbR, \quad \lv F_{i}(x) - F_{i+1}(x) \rv < 1-F_{i+1}(x) \leq \sigma^2 M\rho_{i+1}(x), \quad x \geq \mu_{i+1},     \\
        & \lv F_{j}(x) - F_{j+1}(x) \rv < 1, \quad x \in \mbR, \quad \lv F_{j}(x) - F_{j+1}(x) \rv < F_j(x) \leq \sigma^2 M\rho_{j}(x), \quad x \leq \mu_j,
    \end{aligned}
\end{equation}
we have 
\begin{equation}
    g_{ij}(x) \leq \lbb \begin{aligned}
        & \lv F_{i}(x) - F_{i+1}(x) \rv \frac{\lv F_{j}(x) - F_{j+1}(x) \rv}{\rho_{\theta}(x)} \leq \frac{\lv F_{j}(x) - F_{j+1}(x) \rv}{p_{j}\rho_{j}(x)} = \frac{\sigma^2 M}{p_j},  \quad   x \leq \mu_j,      \\
        & \lv F_{j}(x) - F_{j+1}(x) \rv \frac{\lv F_{i}(x) - F_{i+1}(x) \rv}{\rho_{\theta}(x)} \leq \frac{\lv F_{i}(x) - F_{i+1}(x) \rv}{p_{i+1}\rho_{i+1}(x)} = \frac{\sigma^2 M}{p_{i+1}},\quad   x \geq \mu_{i+1}.
    \end{aligned}\right.
\end{equation}
Notice that $i < j$ implies $\mu_j \geq \mu_{i+1}$, two intervals $(-\infty, \mu_j], [\mu_{i+1}, \infty)$ covers $\mbR$, which provides a global bound on the integrand $g_{ij}(x)$. Decomposing the integrated domain as follows
\begin{equation}
    \begin{aligned}
    \int_{-\infty}^{\infty} g_{ij}(x)dx = & \ \int_{-\infty}^{\mu_{i-1}} g_{ij}(x)dx + \int_{\mu_{i-1}}^{\mu_{j+2}} g_{ij}(x)dx + \int_{\mu_{j+2}}^{\infty} g_{ij}(x)dx \\
    \leq & \ \int_{-\infty}^{\mu_{i-1}} \frac{F_i(x)}{p_{i-1}\rho_{i-1}(x)}dx + \int_{\mu_{j+2}}^{\infty} \frac{1-F_{j+1}(x)}{p_{j+2}\rho_{j+2}(x)}dx + \int_{\mu_{i-1}}^{\mu_{j+2}} g_{ij}(x)dx    \\
    \leq & \ \frac{\sigma^4 M}{dp_{i-1}}e^{-d^2/2\sigma^2} + \frac{\sigma^4 M}{dp_{j+2}}e^{-d^2/2\sigma^2} + \frac{3\sigma^2 dM}{\min\{p_j, p_{i+1} \}},
    \end{aligned}
\end{equation}
where we use the estimate 
\begin{equation}
    \int_{-\infty}^{\mu_{i-1}} \frac{F_i(x)}{p_{i-1}\rho_{i-1}(x)}dx \leq \frac{\sigma^2 M}{p_{i-1}}\int_{-\infty}^{\mu_{i-1}} \frac{e^{-(x-\mu_i)^2/2\sigma^2}}{e^{-(x-\mu_{i-1})^2/2\sigma^2}}dx = \frac{\sigma^4 M}{dp_{i-1}}e^{-d^2/2\sigma^2}.
\end{equation}
\textit{Bound the diagonal terms}:
    Focus on \cref{equ:reduction}, our goal is to reduce the complicated integral over $\mbR$ to a much simpler integral over the interval $[\mu_{i-1}, \mu_i]$. We achieve this goal via three procedures: firstly we restrict the integrated domain to $[\mu_i, \mu_{i+1}]$ as follows
    \begin{equation}\label{equ:step1}
        \begin{aligned}
            & \ \int_{\mbR} \frac{\lp F_i - F_{i + 1} \rp^2}{\rho_{\theta}}dx - \int_{\mu_i}^{\mu_{i+1}} \frac{\lp F_i - F_{i + 1} \rp^2}{\rho_{\theta}}dx  \\
            = & \ \int_{-\infty}^{\mu_{i}} \frac{\lp F_i - F_{i + 1} \rp^2}{\rho_{\theta}}dx + \int_{\mu_{i+1}}^{\infty} \frac{\lp F_i - F_{i + 1} \rp^2}{\rho_{\theta}}dx       \\
            \leq & \ \int_{-\infty}^{\mu_{i}} \frac{\sigma^4M^2\rho_i^2(x)}{p_i \rho_i(x)}dx + \int_{\mu_{i+1}}^{\infty} \frac{\sigma^4M^2\rho_{i+1}^2(x)}{p_{i+1} \rho_{i+1}(x)}dx = \sigma^4M^2 \lp \frac{1}{2p_i} + \frac{1}{2p_{i+1}}\rp.
        \end{aligned} 
    \end{equation}
    Secondly, we simplify the denominator as
    \begin{equation}\label{equ:step2}
        \frac{1}{M_1}\int_{\mu_i}^{\mu_{i+1}} \frac{\lp F_i - F_{i + 1} \rp^2}{p_i\rho_i+p_{i+1}\rho_{i+1}}dx \leq \int_{\mu_i}^{\mu_{i+1}} \frac{\lp F_i - F_{i + 1} \rp^2}{\rho_{\theta}}dx \leq \int_{\mu_i}^{\mu_{i+1}} \frac{\lp F_i - F_{i + 1} \rp^2}{p_i\rho_i+p_{i+1}\rho_{i+1}}dx,
    \end{equation}
    with 
    \begin{equation}
        \begin{aligned}
            M_1 = \max_{x\in[\mu_i, \mu_{i+1}]}\frac{\rho_{\theta}(x)}{p_i\rho_i(x)+p_{i+1}\rho_{i+1}(x)} \leq 1 + \frac{e^{-d^2/2\sigma^2}}{\min\{p_i, p_{i+1}\}}.
        \end{aligned}
    \end{equation}
    Lastly, we deal with the numerator:
    \begin{equation}
        \begin{aligned}
        & \ \int_{\mu_i}^{\mu_{i+1}} \frac{\lp F_i - F_{i + 1} \rp^2}{p_i\rho_i+p_{i+1}\rho_{i+1}}dx - \int_{\mu_i+\sqrt{\sigma}}^{\mu_{i+1}-\sqrt{\sigma}} \frac{\lp F_i - F_{i + 1} \rp^2}{p_i\rho_i+p_{i+1}\rho_{i+1}}dx     \\
       =  & \ \int_{\mu_i}^{\mu_{i}+\sqrt{\sigma}} \frac{\lp F_i - F_{i + 1} \rp^2}{p_i\rho_i+p_{i+1}\rho_{i+1}}dx + \int_{\mu_{i+1}-\sqrt{\sigma}}^{\mu_{i+1}} \frac{\lp F_i - F_{i + 1} \rp^2}{p_i\rho_i+p_{i+1}\rho_{i+1}}dx   \\
       \leq & \ \frac{\sqrt{\sigma}}{\frac{p_i}{\sqrt{2\pi}\sigma}e^{-(\sqrt{\sigma})^2/2\sigma^2}} + \frac{\sqrt{\sigma}}{\frac{p_{i+1}}{\sqrt{2\pi}\sigma}e^{-(\sqrt{\sigma})^2/2\sigma^2}} = \lp\frac{1}{p_i} + \frac{1}{p_{i+1}}\rp\sqrt{2\pi\sigma^3}e^{1/2\sigma}.
        \end{aligned}
    \end{equation}
    In the interval $[\mu_i+\sqrt{\sigma}, \mu_{i+1}-\sqrt{\sigma}]$, we have following estimation:
    \begin{equation}
       1-\sqrt{\frac{2\sigma}{\pi}}e^{-1/2\sigma} \leq (1-\sqrt{\frac{\sigma}{2\pi}}e^{-1/2\sigma})^2 \leq \lp \Phi(\sqrt{\sigma}/\sigma) - \Phi(-\sqrt{\sigma}/\sigma) \rp^2 = \lp F_i(x) - F_{i + 1}(x) \rp^2 \leq 1,
    \end{equation}
    which translates to 
    \begin{equation}\label{equ:step3}
        \lp 1-\sqrt{\frac{2\sigma}{\pi}}e^{-1/2\sigma} \rp \int_{\mu_i+\sqrt{\sigma}}^{\mu_{i+1}-\sqrt{\sigma}} \frac{\lp F_i - F_{i + 1} \rp^2}{p_i\rho_i+p_{i+1}\rho_{i+1}}dx \leq \int_{\mu_i+\sqrt{\sigma}}^{\mu_{i+1}-\sqrt{\sigma}} \frac{\lp F_i - F_{i + 1} \rp^2}{p_i\rho_i+p_{i+1}\rho_{i+1}}dx.
    \end{equation}
    Now, combining \cref{equ:step1}, \cref{equ:step2}, and \cref{equ:step3} together we conclude the bound of the diagonal term.
\end{proof}
\begin{proof}[Proof of \Cref{prop:Delta-2}]
The proof technique is a variant of the Laplace's method in asymptotic analysis. We first prove the Laplace case as follows
\begin{equation}
    \begin{aligned}
       & \ \lim_{\sigma\rightarrow 0} 2\sigma\int_0^d \frac{dx}{p_ie^{-x/\sigma}+\frac{p_{i+1}}{k}e^{-(d-x)/k\sigma}}        \\
       = & \ \lim_{\sigma\rightarrow 0} 2\sigma\int_{-l}^{d-l} \frac{dx}{p_ie^{-(x+l)/\sigma}+\frac{p_{i+1}}{k}e^{-(d-l-x)/k\sigma}}       \\
       = & \ \lim_{\sigma\rightarrow 0} 2\sigma^2\int_{-\infty}^{\infty} \frac{du}{p_ie^{-u-l/\sigma}+\frac{p_{i+1}}{k}e^{u/k-(d-l)/k\sigma}}        \\
       = & \ \lim_{\sigma\rightarrow 0} \frac{2\sigma^2}{p_i}e^{l/\sigma}\int_{-\infty}^{\infty} \frac{du}{e^{-u}+e^{u/k}}        \\
       = & \ \lim_{\sigma\rightarrow 0} \frac{2\sigma^2}{p_i}e^{l/\sigma}\int_{-\infty}^{\infty} \frac{e^udu}{1+e^{u(k+1)/k}}        \\
       = & \ \lim_{\sigma\rightarrow 0} \frac{2\sigma^2}{p_i}e^{l/\sigma}\int_{0}^{\infty} \frac{dy}{1+y^{(k+1)/k}}.        \\
    \end{aligned}
\end{equation}
In the first, the second, and the last deduction, we use the change of variables $x \rightarrow x+l, x \rightarrow u = \frac{x}{\sigma}, u \rightarrow y = e^u$ respectively. The first approximation we used is changing integration domain from $[-l, d-l]$ to $\lp -\infty, \infty \rp$, which can be absorbed in an $O(1)$ factor. Compared with the main part, this is sufficiently small to be ignore.

Next, we prove the asymptotic formula in the Gaussian case following the same derivation as the Laplace case:
\bequn
	\begin{aligned}
		& \lim_{\sigma \rightarrow 0} \sqrt{2\pi}\sigma\int_{0}^{d}\frac{dx}{p_{i} e^{-x^2/2\sigma^2} + \frac{p_{i+1}}{k}e^{-\lp d - x \rp^2/2\lp k\sigma \rp^2}} 		\\
        = & \lim_{\sigma \rightarrow 0} \sqrt{2\pi}\sigma\int_{-l}^{d-l}\frac{dx}{p_{i} e^{-(x+l)^2/2\sigma^2} + \frac{p_{i+1}}{k}e^{-\lp d - l - x \rp^2/2\lp k\sigma \rp^2}} 		\\
        = & \lim_{\sigma \rightarrow 0} \sqrt{2\pi}\sigma^2\int_{-\infty}^{\infty}\frac{du}{p_{i} e^{-(u+l/\sigma)^2/2} + \frac{p_{i+1}}{k}e^{-\lp (d - l)/\sigma - u \rp^2/2k^2}} 		\\
		= & \lim_{\sigma \rightarrow 0} \frac{\sqrt{2\pi}\sigma^2}{p_{i}} e^{\half (l/\sigma)^2}\int_{-\infty}^{\infty} \frac{du}{ e^{ - u^2/2 - ul/\sigma } + e^{u (d - l)/(k^2\sigma) - \lp u/k \rp^2/2  }}	\\
		= & \lim_{\sigma \rightarrow 0} \frac{\sqrt{2\pi}\sigma^2}{p_{i}} e^{\half (l/\sigma)^2}\int_{-\infty}^{\infty} \frac{du}{ e^{- ul/\sigma } + e^{u (d - l)/(k^2\sigma) }}	\\
        = & \lim_{\sigma \rightarrow 0} \frac{\sqrt{2\pi}\sigma^2}{p_{i}} e^{\half (l/\sigma)^2}\int_{-\infty}^{\infty} \frac{e^{ul/\sigma }du}{ 1 + e^{u (d - l)/(k^2\sigma)+ul/\sigma }}	\\
        = & \lim_{\sigma \rightarrow 0} \frac{\sqrt{2\pi}\sigma^3}{p_{i}l} e^{\half (l/\sigma)^2}\int_{0}^{\infty} \frac{dy}{ 1 + y^{1+(d-l)/k^2l} }	\\
        = & \lim_{\sigma \rightarrow 0} \frac{\sqrt{2\pi}\sigma^3}{p_{i}l} e^{\half (l/\sigma)^2}\int_{0}^{\infty} \frac{dy}{ 1 + y^{(k+1)/k} }
	\end{aligned}
\eequn
Most of the derivation follows exactly as the Laplace case. There are two different steps: one is in the fourth derivation where we drop the quadratic term $-u^2/2, \lp u/k \rp^2/2$, which does not depend on the large variable $l/\sigma$. This is a standard technique in reducing asymptotic integral, c.f. \cite{murray2012asymptotic}. Later, we will use perturbation analysis on this term in \Cref{lem:approx} to derive the higher-order scaling Wasserstein metric in \Cref{thm:second-order-metric}. Moreover, in the last deduction, we use the fact that the integral is convergent uniformly over $t$ and the following limit:
\bequn
	\lim_{t \rightarrow \infty}\frac{d - l}{k^2l} = \lim_{\sigma \rightarrow 0}\frac{d - l}{k^2l} = \frac{1}{k},
\eequn
which is derived using \cref{equ:match-cond-prop} to obtain 
\begin{equation}\label{equ:l-gaussian}
    l = \frac{d}{k + 1} + \frac{k\sigma^2}{d}\log \frac{kp_{i}}{p_{i+1}} + o\lp \sigma^2 \rp.
\end{equation} 
This is the second difference with the proof in Laplace mixture model and will also be studied in \Cref{lem:approx} for higher-order metric.

\end{proof}

\section{Generalizations of the scaling Wasserstein information matrices}\label{generalize}

The derivation of the scaling limit in \Cref{limit} requires the unnecessary \Cref{assm} and other constraint. It can only be applied to the 1D case where gaps between adjacent components are the same, which does not hold in most cases. In this section, we study several generalizations of this model in separate subsections.
\subsection{Inhomogeneous Gap}

In this section, we first consider the case where gaps between adjacent means are not all the same, i.e., $\mu_{i + 1} - \mu_i = d_i$s are not constants. We call these models 1D inhomogeneous GMMs. 

We first argue why the method for homogeneous case does not generalize to current setting: if we apply the asymptotic formula in homogeneous case, the factors $K_i\lp \sigma \rp = \sqrt{2\pi^3}\frac{\sigma^3}{d}e^{\half \lp \frac{d_i}{2\sigma}\rp^2}$ appearing in scaling limits of $\lp G_W\lp \theta; \sigma \rp \rp_{i,i}$ is not the same for all components. Consequently, different elements of WIMs have different scaling factors under the limit $\sigma \rightarrow 0$ and there does not exist a scaling factor $K(\sigma)$ such that \cref{equ:W-limit-matrix} holds with a non-trivial limit. Therefore, we consider GMMs with different scaling behaviors of each component's variances, i.e., $\sigma_i$s are not the same but are related by some quantities. We will show that WIMs of models behave well under such a scaling limit. The corresponding scaling limit is stated below.

\begin{Thm}\label{thm:inW-limit}
	For an 1D inhomogeneous GMM with the gaps, variances given by $\mu_{i + 1} - \mu_i = d_i, \sigma_i$ respectively, suppose the following relation holds:
	\bequ\label{equ:match-cond}
		\frac{d_1}{s_1 + s_2}  = \frac{d_2}{s_2 + s_3} = \cdots = \frac{d_{N - 1}}{s_{N - 1} + s_N} = d, \quad \sigma_i = s_i \sigma \quad \forall i \in [N].
	\eequ
	The scaling limit of WIMs is given by
	\bequ\label{equ:inWasserstein-limit}
		\begin{aligned}
		\lim_{\sigma \rightarrow 0} \frac{\lp G_W\lp \theta; \sigma \rp \rp_{ij}}{K_{in}\lp \sigma \rp} &= \delta_{ij}g\lp s_{i+1}/s_i \rp \lp \frac{s_{i+1}}{p_{i+1}}\rp^{\frac{s_{i+1}}{s_{i+1} + s_i}}\lp \frac{s_{i}}{p_i}\rp^{\frac{s_i}{s_{i+1} + s_i}} s_i,
		\end{aligned}\tag{in-W-limit}
	\eequ
    with 
	\bequn
		\begin{aligned}	
  K_{in}\lp \sigma \rp = \frac{\sqrt{2\pi}\sigma^3}{d}e^{\half \lp \frac{d}{\sigma}\rp^2},
		\end{aligned}
	\eequn
\end{Thm}
	It is simple to derive that the inhomogeneous scaling limit \cref{equ:inWasserstein-limit} contains homogeneous scaling limit \cref{equ:W-limit} as a special case. Suppose \Cref{assm} holds, a solution of the matching condition \cref{equ:match-cond} is given by $s_i = k_i = 1, \sigma_i = \sigma, d_i = 2d, \ \forall i \in [N]$, and we have
	\bequn
		\lim_{\sigma \rightarrow 0} \frac{\lp G_W\lp \theta; \sigma \rp \rp_{ii}}{K_{in}\lp \sigma \rp} = g\lp 1 \rp \lp \frac{1}{p_{i+1}}\rp^{\half}\lp \frac{1}{p_i}\rp^{\half} \rightarrow \frac{\pi/2}{\sqrt{p_ip_{i + 1}}}.
	\eequn
    Comparing with \cref{equ:W-limit}, the additional $\pi/2$ factor comes from our different definition of $K(\sigma), K_{in}(\sigma)$. The calculation of the scaling limit requires the full strength of \Cref{prop:Delta-2} to deal with in-homogeneous gaps.

\begin{proof}[Proof of \Cref{thm:inW-limit}]
Using the same argument in the proof of \Cref{thm:W-limit}, we can focus on the diagonal terms of the WIM and consider the integral $\int_0^{d_i} \frac{dx}{p_{i}\rho \lp x; \sigma_i \rp + p_{i+1} \rho \lp d_i - x; \sigma_{i+1} \rp}$. Next, we use the \Cref{prop:Delta-2} with $\sigma=\sigma_i, l = l_i, k = k_i = s_{i+1}/s_i$ to conclude that
\bequ
	\lim_{\sigma \rightarrow 0} \frac{\int_{\mu_i}^{\mu_{i+1}} \frac{dx}{p_{i}\rho(x-\mu_i;\sigma_i) + p_{i+1} \rho(x-\mu_{i+1};\sigma_{i+1})}}{\frac{\sqrt{2\pi} \sigma_i^3 }{p_{i}l_i}e^{\half \lp l_i/\sigma_i \rp^2}} = g \lp k_i \rp,
\eequ
where $l_i$ is given by \cref{equ:l-gaussian}. Using \cref{equ:l-gaussian} and \cref{equ:match-cond}, we obtain
\begin{equation}
    \begin{aligned}
    \frac{\sqrt{2\pi} \sigma_i^3 }{p_{i}l_i}e^{\half \lp l_i/\sigma_i \rp^2} = \lp \frac{\sqrt{2\pi} s_i^2\sigma^3 }{p_{i}d}+ O(\sigma^5)\rp\lp \frac{s_{i+1}p_i}{s_ip_{i+1}}\rp^{\frac{s_{i+1}}{s_i+s_{i+1}}}e^{\half \lp d/\sigma \rp^2}.
    \end{aligned}
\end{equation}
Combining with previous results we conclude the proof.
\end{proof}

We use the following example to illustrate the geometric structure of scaling Wasserstein metric on inhomogeneous GMMs.
\begin{Ex}[2-d inhomogeneous GMM]
	We consider here a 3-component inhomogeneous GMM with $\mu_1 = -1, \mu_2 = 0, \mu_3 = 2$, i.e. $d_1 = 1 \neq 2 = d_2$. If we choose the same variance for all the components, i.e. $\mcN \lp \mu_i, \sigma \rp$. Then, following the discussion of homogeneous models, the scaling WIMs is given by
\bequn
	\begin{aligned}
		\lim_{\sigma \rightarrow 0}\frac{\lp G_W \lp \theta; \sigma\rp\rp_{11}}{\sqrt{2\pi^3}\sigma^3/d_1  e^{\half \lp \frac{d_1}{2\sigma}\rp^2}} & = \frac{1}{\sqrt{p_1p_2}}, 			\\
		\lim_{\sigma \rightarrow 0}\frac{\lp G_W \lp \theta; \sigma\rp\rp_{22}}{\sqrt{2\pi^3}\sigma^3/d_2  e^{\half \lp \frac{d_2}{2\sigma}\rp^2}} & = \frac{1}{\sqrt{p_2p_3}}.
	\end{aligned}
\eequn
However, plugging $d_1 = 1, d_2 = 2$ in the scaling limits of $\lp G_W \lp \theta; \sigma\rp\rp_{11}, \lp G_W \lp \theta; \sigma\rp\rp_{22}$ means they diverge with different speed $e^{\half \lp \frac{d_1}{2\sigma}\rp^2}, e^{\half \lp \frac{d_2}{2\sigma}\rp^2}$ as $\sigma \rightarrow 0$. Consequently, we can not find a common scaling factor $K\lp \sigma \rp$ to normalize the metric in the homogeneous sense, c.f. \cref{equ:scaling-limit}. Therefore, we need to choose different variance scaling so that a common scaling factor exists. To start with, we choose variances of both the first and the second components to be the same, namely
\bequn
	\sigma_1 = \sigma_2 = \sigma.
\eequn
Then we know that the matrix element $\lp G_W \lp \theta; \sigma\rp\rp_{11}$ has the following scaling limit 
\bequn
	\lim_{\sigma \rightarrow 0}\frac{\lp G_W \lp \theta; \sigma\rp\rp_{11}}{\sqrt{2\pi^3}\sigma^3  \frac{e^{\half \lp \frac{d_1}{2\sigma}\rp^2}}{d_1}} = \frac{1}{\sqrt{p_1p_2}}.
\eequn
Next, we need to choose a specific $\sigma_3$ such that the matrix element $\lp G_W \lp \theta; \sigma\rp\rp_{22}$ have the same scaling factor. Using the conclusion for the inhomogeneous models \Cref{thm:inW-limit}, we conclude that 
\bequn
	\lim_{\sigma \rightarrow 0} \frac{\lp G_W \lp \theta; \sigma\rp\rp_{22}}{\frac{\sqrt{2\pi}\sigma^3 \lp \frac{\sigma_3 + \sigma}{\sigma} \rp}{d_2} e^{\half \lp \frac{d_2}{\sigma + \sigma_3}  \rp^2}} =  \frac{\lp \frac{\sigma_3}{\sigma} \rp^{\frac{\sigma_3}{\sigma + \sigma_3}}}{p_{2}^{\frac{\sigma}{\sigma + \sigma_3}}p_3^{\frac{\sigma_3}{\sigma + \sigma_3}}}\lp g\lp \frac{\sigma_3}{\sigma} \rp + \frac{\sigma_3}{\sigma}f \lp \frac{\sigma_3}{\sigma} \rp \rp.
\eequn
Thus we require
\bequn
	\frac{d_2}{\sigma + \sigma_3} = \frac{d_1}{2\sigma},
\eequn
which is exactly the matching condition \cref{equ:match-cond} we stated before. With $d_2 = 2d_1$, we have 
\bequn
	\sigma_3 = 3\sigma.
\eequn
The scaling limit of the metric inner product can be further given as
\bequn
	\begin{aligned}
		\lim_{\sigma \rightarrow 0}\frac{\lp G_W \lp \theta; \sigma\rp\rp_{11}}{2\sqrt{2\pi}\sigma^3  e^{1/8\sigma^2}} & = \frac{\pi}{2\sqrt{p_1p_2}}, 			\\
		\lim_{\sigma \rightarrow 0} \frac{\lp G_W \lp \theta; \sigma\rp\rp_{22}}{2\sqrt{2\pi}\sigma^3 e^{1/8\sigma^2}} & = \frac{3^{\frac{3}{4}}}{p_{2}^{\frac{1}{4}}p_3^{\frac{3}{4}}}g(3).
	\end{aligned}
\eequn
Combining all the result in hand, we have that the limit metric for this GMM is given by
\bequn
	G_{\wtd W}^{(in)} \lp \theta \rp = \lim_{\sigma \rightarrow 0}\frac{G_W \lp \theta; \sigma\rp}{2\sqrt{2\pi}\sigma^3  e^{1/8\sigma^2}} = \begin{pmatrix}
		\frac{\pi}{2\sqrt{p_1p_2}} & 0		\\
		0 & \frac{3^{\frac{3}{4}}}{p_{2}^{\frac{1}{4}}p_3^{\frac{3}{4}}}g(3)
	\end{pmatrix}.
\eequn
\end{Ex}

\subsection{A Second-order Metric}
In deriving the scaling limit \cref{equ:W-limit}, we take the limit under $\sigma \rightarrow 0$  and analyze the leading order term of the integral which defines the metric. In this subsection, we present a more accurate analysis which also address the second order term in the metric integral. We use this information to form the second extended model, i.e. the second order metric.

More intuitively speaking, we expand the metric $G\lp \theta, \sigma \rp$ w.r.t. $\sigma$, i.e.
\bequn
	\begin{aligned}
		G\lp \theta, \sigma \rp  = & \  G\lp \theta \rp + \frac{\p G\lp \theta, \sigma \rp}{\p \sigma} g\lp \sigma \rp + o\lp g\lp \sigma \rp \rp.
	\end{aligned}
\eequn
\Cref{thm:W-limit} in the first section can be viewed as calculating $\lim_{\sigma \rightarrow 0}G\lp \theta, \sigma \rp = G\lp \theta \rp$. In this section, we try to figure out $\frac{\p G\lp \theta, \sigma \rp}{\p \sigma}$, which should be formally understood as some ``derivative'' of the $G$ w.r.t. $\sigma$. $g\lp \sigma \rp$ is the order of the first expansion term. Notice the above formula is merely an illustration of the main idea instead of the rigor formulation which is deduced below. As before, we first state the result here.
\begin{Thm}\label{thm:second-order-metric}
With the same homogeneous setting and scaling factor $K(\sigma) = \sqrt{2\pi^3}\frac{\sigma^3}{d}e^{\frac{d^2}{8\sigma^2}}$ as \Cref{thm:W-limit}, the scaling WIM of a 1D GMM has following expansion 
\bequ\label{equ:2-metric}
	\begin{aligned}
	\frac{G_{\wtd W}^{(2)}\lp \theta; \sigma \rp}{K\lp \sigma \rp} 
	= & \ \begin{pmatrix}
		\frac{1 + \lp\frac{\pi^2}{2} + \frac{4}{\pi}\log\frac{p_1}{p_{2}}g'(1) + \frac{2}{\pi} \log^2\frac{p_1}{p_{2}}\rp\frac{\sigma^2}{d^2}}{\sqrt{p_1p_2}}  & \cdots & 0		\\
		0 & \cdots & 0	\\
		\vdots & \ddots & \vdots		\\
		 0 & \cdots & \frac{1 + \lp\frac{\pi^2}{2} + \frac{4}{\pi}\log\frac{p_{N-1}}{p_{N}}g'(1) + \frac{2}{\pi} \log^2\frac{p_{N-1}}{p_{N}}\rp\frac{\sigma^2}{d^2}}{\sqrt{p_{N - 1}p_N}} 
	\end{pmatrix}			\\
	= & \ G_{\wtd W}\lp \theta \rp + \frac{\sigma^2}{d^2}\begin{pmatrix}
		\frac{\frac{\pi^2}{2} + \frac{4}{\pi}\log\frac{p_1}{p_{2}}g'(1) + \frac{2}{\pi} \log^2\frac{p_1}{p_{2}}}{\sqrt{p_1p_{2}}} & \cdots & 0		\\
		0 & \cdots & 0	\\
		\vdots & \ddots & \vdots		\\
		 0 & \cdots & \frac{\frac{\pi^2}{2} + \frac{4}{\pi}\log\frac{p_{N-1}}{p_{N}}g'(1) + \frac{2}{\pi} \log^2\frac{p_{N-1}}{p_{N}}}{\sqrt{p_{N - 1}p_N}}
	\end{pmatrix},
	\end{aligned}\tag{2-W-limit}
\eequ
where higher order term $O(\sigma^4)$ is omitted.
\end{Thm}
One should understand this new result as a refinement of the old one \cref{equ:W-limit-matrix}. A similarity shared by this higher-order Wasserstein metric \cref{equ:2-metric} and first-order metric \cref{equ:W-limit-matrix} is that both of them are diagonal. Consequently, inverting the metric $G_W^{-1}$ tensor is straight-forward, making it easy to obtain the gradient flow equations in these models, which will be discussed in \Cref{gradient-flow}. Moreover, the appearance of the fraction $\frac{\sigma}{d}$ is natural as it describes the relative size of the Gaussian components w.r.t. the underlying grid space.

To begin with, we first develop a key estimation frequently used in the derivation below.
\begin{Lem}\label{lem:approx}
	We have following perturbation expansion
	\bequn	
		\begin{aligned}
			& \ \int_{-\infty}^{\infty} \frac{du}{ e^{ - u^2/2 - tu } + e^{tu/k_1 - \lp u/k_2 \rp^2/2  }}	\\
		= & \ \int_{-\infty}^{\infty} \frac{du}{ e^{- ut } + e^{tu/k_1 }} + \frac{1+1/k_2^2}{4t^3}g_2(k_1) + O\lp\frac{1}{t^4}\rp,
		\end{aligned}
	\eequn
 where $g_2(k_1) = \int_{0}^{\infty} \frac{\log^2 v }{ 1 + v^{(k_1+1)/k_1 }}dv$.
\end{Lem}
\begin{proof}
We first expand the ratio of two denominators according to $u$:
\begin{equation}
    \begin{aligned}
    & \ \frac{e^{ - u^2/2 - tu } + e^{tu/k_1 - \lp u/k_2 \rp^2/2}}{e^{- ut } + e^{tu/k_1 }}         \\ 
    = & \ \frac{e^{-tu}\lb 1-\half u^2 + O(u^4)\rb + e^{tu/k_1}\lb 1-\half \lp u/k_2 \rp^2+O(u^4)\rb}{e^{- ut } + e^{tu/k_1 }}    \\
    = & \ 1 + \frac{-e^{-tu}(\half u^2 + O(u^4)) - e^{tu/k_1}(\half \lp u/k_2 \rp^2)+O(u^4)}{1-ut+\half (ut)^2 + 1 + tu/k_1 + \half (tu/k_1)^2 + O(u^3)}       \\
    = & \ 1 - \frac{1+1/k_2^2}{4}u^2 + O(u^3).
    \end{aligned}
\end{equation}
Consequently, we have
\begin{equation}
    \begin{aligned}
			& \ \int_{-\infty}^{\infty} \frac{du}{ e^{ - u^2/2 - tu } + e^{tu/k_1 - \lp u/k_2 \rp^2/2  }}	\\
		= & \ \int_{-\infty}^{\infty} \frac{1}{1 - \frac{1+1/k_2^2}{4}u^2 + O(u^3)}\frac{du}{ e^{- ut } + e^{tu/k_1 }}        \\
        = & \ \int_{-\infty}^{\infty} \frac{du}{ e^{- ut } + e^{tu/k_1 }} + \int_{-\infty}^{\infty} \frac{\frac{1+1/k_2^2}{4}u^2 + O(u^3)}{ e^{- ut } + e^{tu/k_1 }} du       \\
        = & \ \int_{-\infty}^{\infty} \frac{du}{ e^{- ut } + e^{tu/k_1 }} + \frac{1+1/k_2^2}{4t^3}g_2(k_1) + O(\frac{1}{t^4}),       \\
		\end{aligned}
\end{equation}
where the last derivation is based on following calculation
\begin{equation}
    \int_{-\infty}^{\infty} \frac{u^m}{ e^{- ut } + e^{tu/k_1 }}du = \int_{-\infty}^{\infty} \frac{u^m e^{tu}}{ 1 + e^{tu(k_1+1)/k_1 }}du = \frac{1}{t^{m+1}}\int_{0}^{\infty} \frac{\log^m v }{ 1 + v^{(k_1+1)/k_1 }}dv.
\end{equation}
\end{proof}
\begin{proof}[Proof of \Cref{thm:second-order-metric}]
    The second order expansion consists of three parts we will calculate separately. Firstly, we apply \Cref{lem:approx} to the homogeneous case where $k = 1, t = \frac{l}{\sigma}, k_1 = \frac{l}{d-l} = 1+O(\sigma^2), k_2 = k = 1$ to obtain
    \bequn	
		\begin{aligned}
			& \ \int_{-\infty}^{\infty} \frac{du}{ e^{ - u^2/2 - ul/\sigma } + e^{u (d - l)/(k^2\sigma) - \lp u/k \rp^2/2  }}	\\
		= & \ \int_{-\infty}^{\infty} \frac{du}{ e^{- ul/\sigma } + e^{u (d - l)/(k^2\sigma) }} + \frac{\sigma^3}{2l^3}\frac{\pi^3}{8} + O(\sigma^4),
		\end{aligned}
	\eequn
 where the integral of $g_2(1)$ can be formed explicitly to obtain $\frac{\pi^3}{8}$ \cite{191765}. The first term is the main part which integrates to $\frac{\sigma\pi}{2l}$ in the proof of \cref{equ:W-limit}. 
 
 The second approximation appears when we derive
\begin{equation}
    \int_{0}^{\infty} \frac{dy}{ 1 + y^{1+(d-l)/k^2l} }	
        \rightarrow \int_{0}^{\infty} \frac{dy}{ 1 + y^{(k+1)/k} }.
\end{equation}
Recall the expansion of $l$ over $\sigma$ is given by \cref{equ:l-gaussian}, hence we have
\begin{equation}
    \begin{aligned}
    & \ \int_{0}^{\infty} \frac{dy}{ 1 + y^{1+(d-l)/k^2l} }	
        - \int_{0}^{\infty} \frac{dy}{ 1 + y^{(k+1)/k} } \\
        = & \ g(k+\frac{k(k+1)^2\sigma^2}{d^2}\log\frac{kp_i}{p_{i+1}}) - g(k) \\
        = & \ g(1+\frac{4\sigma^2}{d^2}\log\frac{p_i}{p_{i+1}}) - g(1) = \frac{4\sigma^2}{d^2}\log\frac{p_i}{p_{i+1}}\frac{dg}{dx}\Big|_{x=1} + O(\sigma^4).
    \end{aligned}
\end{equation}
This finishes part of the second order expansion of $I_i$, i.e.
\begin{equation}
    I_i = \frac{\sqrt{2\pi}\sigma^2}{p_{i}} \frac{\sigma}{l}e^{\half (l/\sigma)^2}\lp \frac{\pi}{2} + \frac{\sigma^2\pi^3}{4d^2} + \frac{4\sigma^2}{d^2}\log\frac{p_i}{p_{i+1}}\frac{dg}{dx}\Big|_{x=1} + O(\sigma^4)\rp.
\end{equation}
Now, we work on the expansion for the factor contains $l$ as it can be expanded as a power series of $\sigma$. In the homogeneous scenario, \cref{equ:l-gaussian} reduces to exact formula
\begin{equation}
    l = \frac{d}{2} + \frac{\sigma^2}{d}\log \frac{p_{i}}{p_{i+1}}.
\end{equation}
Plug in, one obtains
\begin{equation}
    \frac{\sigma}{l}e^{\half (l/\sigma)^2} = \frac{2\sigma}{d}\sqrt{\frac{p_i}{p_{i+1}}} e^{d^2/8\sigma^2+\lp\frac{\sigma}{d}\log \frac{p_i}{p_{i+1}}\rp^2/2}\lp 1 - \frac{2\sigma^2}{d^2}\log \frac{p_i}{p_{i+1}} + O(\sigma^4)\rp
\end{equation}
Consequently, we have 
 \bequ
	I_i = \sqrt{2\pi^3}\frac{\sigma^3}{d}e^{d^2/8\sigma^2}\frac{1}{\sqrt{p_{i}p_{i+1}}} \lp 1 + \lp\frac{\pi^2}{2} + \frac{4}{\pi}\log\frac{p_i}{p_{i+1}}\frac{dg}{dx}\Big|_{x=1} + \frac{2}{\pi} \log^2\frac{p_i}{p_{i+1}}\rp\frac{\sigma^2}{d^2} + O(\sigma^4)\rp.
\eequ
\end{proof}
Recall in the proof of \Cref{prop:Delta-2} we compare the derivations in Gaussian and Laplace case. The approximation in \Cref{lem:approx} does not exist in Laplace case as its exponent is linear in $x$. Similarly, the approximation related to $\frac{dg}{dx}$ is also not present. Lastly, as the parameter $l$ can be derived in closed-form in Laplace case comparing to the expansion of $l$ in $\sigma$ in Gaussian case, c.f. \cref{equ:l-laplace}, \cref{equ:l-gaussian}, the scaling factor of Laplace mixture model is also exact. Consequently, we conclude that the scaling WIM is accurate to any order of $\sigma$, i.e. all the higher-order expansions vanish. The difference between Laplace and Gaussian cases also indicates that the second-order metric is not universal for all the mixture models. Meanwhile, the same technique is applicable to obtain the metric of any order accuracy.

\subsection{An extended 1D GMM}\label{sec:ext-GMM}
In this section, as promised, we consider an extended GMM which is more akin to GMM in statistical community \cite{reynolds2009gaussian} in the sense that mean parameters are also allowed to vary. Recall in \cref{equ:GMM} each component $\rho_i(x) = \frac{1}{\sqrt{2\pi}\sigma_i}e^{-(x-\mu_i)^2/2\sigma_i^2}$ is fixed a priori and only $\theta$-parameters can vary in the model. In this extended model, the mean variables $\mu_i$s are also included in the model's parameters, thus allowing each component to move along the real axis. We call this ``1D extended GMMs''. Considering its tangent space over $(\theta, \mu)$, it has tangent vectors $\frac{\p }{\p \mu_i}$s associated with each mean parameters. Thus, on this GMM, a basis of tangent space is given by $\lbb \frac{\p}{\p \mu_i}, i \in [N], \frac{\p}{\p \theta_j}, j \in [N-1] \rbb$, which contains exactly a basis of original GMMs, i.e. $\lbb \frac{\p}{\p \theta_j}, j \in [N-1] \rbb$. We will analyze the WIM associated with this larger basis. The results is summarized as follows.
\begin{Thm}\label{thm:extW-limit}
	The WIM in 1D extended homogeneous GMM is given in following block form
	\bequ\label{equ:ex-W-metric}
		\begin{aligned}
			& \lim_{\sigma \rightarrow 0}G_{\wtd W}^{(ext)}\lp \theta, \mu; \sigma \rp = \begin{pmatrix}
				\lp G_{\wtd W}^{(ext)} \rp_{\theta\theta} & \lp G_{\wtd W}^{(ext)} \rp_{\theta\mu}		\\
				\lp G_{\wtd W}^{(ext)} \rp_{\mu\theta} & \lp G_{\wtd W}^{(ext)} \rp_{\mu\mu}
			\end{pmatrix},						\\
			& \lp G_{\wtd W}^{(ext)} \rp_{\theta\theta} = K\lp \sigma \rp \begin{pmatrix}
				\frac{1}{\sqrt{p_1p_2}} & 0 & \cdots & 0		\\
				0 & \frac{1}{\sqrt{p_2p_3}} & \cdots & 0		\\
				\vdots & \vdots & \ddots & \vdots \\
				0 & 0 & \cdots & \frac{1}{\sqrt{p_{N - 1}p_N}}
			\end{pmatrix}, \quad \lp G_{\wtd W}^{(ext)} \rp_{\mu\mu} = \begin{pmatrix}
				p_{1} & 0 & \cdots & 0		\\
				0 & p_{2} & \cdots & 0		\\
				\vdots & \vdots & \ddots & \vdots \\
				0 & 0 & \cdots & p_{N}
			\end{pmatrix},							\\
			& \lp \lp G_{\wtd W}^{(ext)} \rp_{\mu\theta} \rp^T = \lp G_{\wtd W}^{(ext)} \rp_{\theta\mu} = \begin{pmatrix}
				\frac{\mu_2 - \mu_1}{2} & \frac{\mu_2 - \mu_1}{2} & 0 & \cdots & 0 & 0		\\
				0 & \frac{\mu_3 - \mu_2}{2} & \frac{\mu_3 - \mu_2}{2} & \cdots & 0 & 0		\\
				\vdots & \vdots & \vdots & \ddots & \vdots & \vdots \\
				0 & 0 & 0 & \cdots & \frac{\mu_{N - 1} - \mu_{N - 2}}{2} & 0		\\
				0 & 0 & 0 & \cdots & \frac{\mu_N - \mu_{N - 1}}{2} & \frac{\mu_N - \mu_{N - 1}}{2}
			\end{pmatrix}.			\\
		\end{aligned}
	\eequ
\end{Thm}
	Let us take a close look at each block. The left upper block $\lp G_{\wtd W}^{(ext)} \rp_{\theta\theta}$ is significantly greater than the other three blocks as $\sigma \rightarrow 0$. In other words, the WIM of this model has a ``multiscale'' feature. Simply dividing the whole matrix by the factor $K(\sigma)$ will result in a degenerate matrix. If we want to obtain a nondegenerate scaling limit of WIM, it suffices to rescale the tangent vectors $\frac{\p}{\p \theta_j}$ by factor $\frac{1}{\sqrt{K\lp \sigma \rp}}$ while remain $\frac{\p}{\p \mu_i}$ unchanged. Consequently, the ``rescaled'' WIM for the extended GMM with tangent vectors given by $\frac{1}{\sqrt{K\lp \sigma \rp}}\frac{\p}{\p \theta_i}, i = 1,2, \cdots, N - 1$ and $\frac{\p}{\p \mu_i}$ reads
	\bequn
		\begin{pmatrix}
			\begin{pmatrix}
				\frac{1}{\sqrt{p_1p_2}} & 0 & \cdots & 0		\\
				0 & \frac{1}{\sqrt{p_2p_3}} & \cdots & 0		\\
				\vdots & \vdots & \ddots & \vdots \\
				0 & 0 & \cdots & \frac{1}{\sqrt{p_{N - 1}p_N}}
			\end{pmatrix} & \mathbf{0}			\\
			\mathbf{0} & \begin{pmatrix}
				p_{1} & 0 & \cdots & 0		\\
				0 & p_{2} & \cdots & 0		\\
				\vdots & \vdots & \ddots & \vdots \\
				0 & 0 & \cdots & p_{N}
			\end{pmatrix}
		\end{pmatrix}.
	\eequn
	Notice the rescaled WIM is also a diagonal matrix, consistent with all the previous scaling limit. The off-diagonal terms vanish due to the scaling factor. Before the proof of \Cref{thm:extW-limit}, we state a relation between two subsets of basis $\lbb \frac{\p }{\p \mu_i} \rbb, \lbb \frac{\p }{\p p_j} \rbb$ in term of the Fisher-Rao geometry and Wasserstein geometry. This result is used for the derivation of the scaling limit of WIM in extended GMM.\begin{Thm}\label{thm:WIG1}
	For the extended GMM, we have the following relationship between the WIM of the set of tangent vectors $\frac{\p}{\p \mu_i}$s and Fisher information matrix of the set of tangent vectors $\frac{\p}{\p \theta_i}$s
	\bequ\label{WIG-eq1}
		G_F = \Sigma G_W \Sigma^T,
	\eequ
	where $\lp G_F \rp_{ij} = g_F\lp \frac{\p}{\p \theta_i}, \frac{\p}{\p \theta_j} \rp$, $\lp G_W \rp_{ij} = g_W\lp \frac{\p}{\p \mu_i}, \frac{\p}{\p \mu_j} \rp$ and the matrix $\Sigma \in \mbR^{N-1,N}$ appears above is given by
	\bequn
		\Sigma = \begin{pmatrix}
			- \frac{1}{p_1} & \frac{1}{p_2} & 0 & \cdots & 0 & 0		\\
			0 & - \frac{1}{p_2} & \frac{1}{p_2} & \cdots & 0	& 0		\\
			\vdots & \vdots & \ddots & \cdots	& \cdots & \cdots				\\
			0 & 0 & \cdots & - \frac{1}{p_{N - 2}} & \frac{1}{p_{N - 1}} & 0		\\
			0 & 0 & 0 & \cdots & - \frac{1}{p_{N - 1}} & \frac{1}{p_N}
		\end{pmatrix}.
	\eequn 
\end{Thm}
\begin{proof}

According to \Cref{thm:WIG1}, we have that the Wasserstein score functions associated with this basis are given by
\bequn
	\begin{aligned}
		\nabla \Phi_{\theta_i}^W(x) & = \frac{F_{i}(x) - F_{i + 1}(x)}{\rho(x)}, \quad \forall i = 1,2,...,N-1,    \\
		\nabla \Phi_{\mu_i}^W(x) & = p_i \frac{\rho_i(x) }{\rho(x)}, \quad \forall i = 1,2,...,N.
	\end{aligned}
\eequn
We give the metric tensor $G_W$ in a block form, whose left-upper part has been shown in ordinary GMM computation, i.e. \Cref{thm:W-limit}. For the right-lower part, the proof is the same
\bequn
	\begin{aligned}
		& \ \lim_{\sigma \rightarrow 0}g_W\lp \p_{\mu_i}, \p_{\mu_j} \rp		\\
		= & \ \lim_{\sigma \rightarrow 0}p_i p_j\int \frac{ \rho_i(x)\rho_j(x)}{\rho(x)}dx			\\
		= & \ \lim_{\sigma \rightarrow 0}p_j\int_{\lp \mu_i - \epsilon, \mu_i + \epsilon \rp} \rho_j(x) dx			\\
		= & \ \lbb \begin{aligned}
			& 0, \quad i \neq j,			\\
			& p_i, \quad i = j,
		\end{aligned}\rpt
	\end{aligned}
\eequn
where the $\epsilon$ appears in the formula is a constant depends on $\sigma$ s.t. $\lim_{\sigma \rightarrow 0}\epsilon = 0$. Notice that the derivation above is similar to the scaling limit of the Fisher-Rao metric. This is attribute to the linear relation between the Wasserstein geometry and Fisher-Rao geometry we posed in the \Cref{thm:WIG1}.

The last part is to compute the Wasserstein inner product between these two basis sets. With the help of the closed-forms of these two Wasserstein score functions, we conclude that $\lim_{\sigma \rightarrow 0}\mbE [\nabla \Phi_{\mu_i}^W \cdot \nabla \Phi_{\theta_j}^W] = 0$ with $i < j, j < i - 1$. The proof is simply a common argument of the decomposition trick stated in \cite{li2023wasserstein}. For the other term, we have
\bequn
	\begin{aligned}
		& \ \lim_{\sigma \rightarrow 0}g_W\lp \p_{\mu_i}, \p_{\theta_i} \rp 		\\
		= & \ \lim_{\sigma \rightarrow 0}p_i\int \frac{ \rho_i(x)}{\rho(x)} \lp F_{i}(x) - F_{i + 1}(x) \rp	 dx	\\
		= & \ \lim_{\sigma \rightarrow 0}p_i\int_{\mu_i}^{\mu_{i + 1}} \frac{ \rho_i(x)}{\rho} dx		\\
		= & \ \frac{\mu_{i + 1} - \mu_i}{2},
	\end{aligned}
\eequn
where the second equality holds by the weak convergence of Gaussian $\mcN\lp \mu, \sigma \rp \rightarrow \delta_{\mu}$ as $\sigma \rightarrow 0$. The last equality holds by the decomposition of the integral and estimation on each part of the integral
\bequn
	\begin{aligned}
		& \lim_{\sigma \rightarrow 0}p_i\int_{\mu_i}^{\frac{\mu_{i + 1} + \mu_i}{2}} \frac{ \rho_i(x)}{\rho(x)} dx \rightarrow \lim_{\sigma \rightarrow 0}\int_{\mu_i}^{\frac{\mu_{i + 1} + \mu_i}{2}} dx = \frac{\mu_{i + 1} - \mu_i}{2},		\\
		& \lim_{\sigma \rightarrow 0}p_i\int_{\frac{\mu_{i + 1} + \mu_i}{2}}^{\mu_{i + 1}} \frac{ \rho_i(x)}{\rho(x)} dx \rightarrow \lim_{\sigma \rightarrow 0}\int_{\mu_i}^{\frac{\mu_{i + 1} + \mu_i}{2}} 0 \cdot dx = 0.
	\end{aligned}
\eequn
\end{proof}
\begin{proof}[Proof of \Cref{thm:WIG1}]
	
	Using the language of our previous paper, the Fisher score functions correspond to the tangent vectors $\frac{\p}{\p \theta_i}$ given by
	\bequn
		\Phi_F^i(x) = \frac{\rho_{i + 1}(x) - \rho_i(x)}{\rho(x)}, \quad \forall i = 1,2,...,N-1.
	\eequn 
	
	While, in term of Wasserstein score function for $\frac{\p}{\p \mu_i}$, we have
	\bequn
		\nabla \Phi_W^i(x) = p_i\frac{\rho_i(x)}{\rho(x)}, \quad \forall i = 1,2,...,N,
	\eequn
	are associated with the tangent vector $\frac{\p}{\p \mu_i}$. Consequently, we conclude that the Fisher score function and Wasserstein score function can be connected via the following linear relation
	\bequn
		\Phi_F^i = \Sigma \nabla \Phi_W^i.
	\eequn
	Thus, by the dual formulation of the statistical metric tensor \cite{li2023wasserstein}, namely
	\bequn
		G_F = \mbE \lb \Phi_F \cdot \Phi_F^T \rb, \quad G_W = \mbE \lb \nabla \Phi_W \cdot \nabla \Phi_W^T \rb,
	\eequn
	we conclude that the linear relation between the FIM and WIM holds.
\end{proof}
\begin{Rem}
The above theorem points out the fact that on each point $\theta$ in this parametric statistical model, we have two disjoint tangent subspace $V_{\theta}^F$ and $V_{\theta}^W$, $V_{\theta}^F \cap V_{\theta}^W = \{0\}$. The WIM $G_W$ on $V_{\theta}^W$ has a linear relation with the Fisher information matrix on $V_{\theta}^F$.
\end{Rem}

\section{Gradient flows on various scaling Wasserstein geometry}\label{gradient-flow}

Gradient flows have been intensively studied in the community of optimal transport for several decades \cite{JKO, MAAS20112250}. In this section, we study the gradient flows under the scaling Wasserstein metric. We derive explicit forms of gradient flows under this scaling geometry and consider its connection with gradient flows on density manifold. We first give a brief summary of gradient flow on density manifold, see detailed studies in \cite{Ambrosio2008Gradient, otto2000generalization}.
\subsection{Review of gradient flows on density manifold}\label{general-gf}
First, we provide the general formulation of Wasserstein gradient flow on density manifold. Given a functional on the density manifold $\mcE: \mcP \rightarrow \mbR$, its Wasserstein gradient is given by \cite{villani2021topics}
\bequ
		\nabla_W F(\rho) = -\nabla \cdot \lp \rho \nabla \lp \frac{\delta \mcE}{\delta \rho}\rp\rp.
\eequ
Consequently, the gradient flow of a potential functional $\mcV(\rho) = \mbE_{X\sim\rho}V(X)$ is given by
\bequ\label{potential-flow}
    \p_t\rho(x) = -\nabla_W \mcV(\rho) = \nabla \cdot \lp \rho \nabla \lp \frac{\delta \mcV}{\delta \rho}\rp\rp = \nabla \cdot (\rho \nabla V(x)).
\eequ
Moreover, the gradient flow of the negative entropy functional $\mcH(\rho) = \mbE_{X\sim\rho}\log \rho(X)$ is given by the heat equation
\bequ
    \p_t\rho(x) = \nabla \cdot \lp \rho \nabla \lp \frac{\delta \mcH}{\delta \rho}\rp\rp = \nabla \cdot (\rho \nabla (\log\rho(x) + 1))) = \Delta \rho.
\eequ
In the next subsection, we will derive parametric gradient flow equation under the scaling Wasserstein metric.
\subsection{Gradient flows on scaling Wasserstein metric}

In this section, we derive analytic forms of gradient flows in different scaling Wasserstein geometries introduced before. We first state the equation of gradient flow of a general functional $F$ on a parametric model $\Theta \subset \mbR^{N - 1}$.
\begin{Prop}
Consider the parametric model $\Theta \subset \mbR^{N - 1}$ with the scaling Wasserstein metric $G_{\wtd W}\lp \theta \rp$ \cref{equ:scaling-limit} defined on it. Given a function on this space $F: \Theta \rightarrow \mbR$, the gradient flow of this function is given by
\bequ\label{equ:gradient-flow-eq}
	\begin{aligned}
		\dot{\theta} & \ = - G_{\wtd W}^{-1}\lp \theta \rp \nabla_{\theta} V\lp \theta \rp,			\\
		\dot{\theta_i} & \ = - \sqrt{p_ip_{i + 1}} \lp \nabla_{\theta} V\lp \theta \rp \rp_i,		\\
		\dot{p}_i & \ = - \sqrt{p_ip_{i - 1}} \lp \nabla_{\theta} V\lp \theta \rp \rp_{i - 1} + \sqrt{p_ip_{i + 1}} \lp \nabla_{\theta} V\lp \theta \rp \rp_i.
	\end{aligned}
\eequ
\end{Prop}
\begin{proof}
The first equation is exactly the definition of gradient flows in Riemannian manifold. The second equation is a scalar version of the first one and $\lp \nabla_{\theta} V\lp \theta \rp \rp_i$ refers to its $i$-th component. This formula is simplified due to the fact that scaling WIMs are diagonal in homogeneous case, c.f. \Cref{thm:W-limit}. 

However, parameters $\theta_i$s do not have a probabilistic meaning compared to $p_i$s, we rewrite the above gradient flow equation in term of parameters $p_i$s, which are illustrated as ratios of components
\bequn
	\begin{aligned}
		\dot{p}_i & \ =  \dot{\theta}_{i - 1}	- \dot{\theta}_{i}		\\
		& \ = - \sqrt{p_ip_{i - 1}} \lp \nabla_{\theta} V\lp \theta \rp \rp_{i - 1} + \sqrt{p_ip_{i + 1}} \lp \nabla_{\theta} V\lp \theta \rp \rp_i.
	\end{aligned}
\eequn
\end{proof}
Notice that if we divided the metric tensor $G_W$ by this factor $K(\sigma)$, the solution of the gradient flow is changed by a time reparameterization, namely, suppose the solution of 
\bequn
	\dot{\theta} = -G_{\wtd W}\lp \theta \rp^{-1} \nabla_{\theta}H,
\eequn 
is given by $\theta\lp t \rp$. Then, the solution of the scaling gradient flow
\bequn
	\dot{\theta} = - G_{\wtd W}\lp \theta \rp^{-1}K\lp \sigma \rp \nabla_{\theta}H,
\eequn
is exactly $\wtd \theta\lp t \rp = \theta\lp K\lp \sigma \rp t \rp$. Consequently, we have the freedom to discard this scaling factor when studying the properties of the gradient flows.

Using this we can write out the analytic form of arbitrary gradient flows in the probability simplex under scaling Wasserstein metric. In \cite{villani2021topics}, the author introduced three classes of important energy, i.e. internal $\mcU$, potential $\mcV$, and interaction energy $\mcW$ on the density manifold and carefully derive their corresponding gradient flow equations. Hereby we follow the same spirit of \cite{villani2021topics} to list the analytic form of three classes of energies and their gradient flows. To achieve this, we have to calculate the corresponding energy forms on the probability simplex and then apply \cref{equ:gradient-flow-eq}. We take the potential functional as an example: viewing the point on probability simplex as a mixture of Dirac measure, i.e. $\{p_i, i\in [N]\} \rightarrow \sum_{i=1}^N p_i \delta_{\mu_i}$, we can calculate the value of potential functional on this point as
\begin{equation}
    \mcV(\{p_i\}) = \mcV(\sum_{i=1}^N p_i \delta_{\mu_i}) = \int V(x)\sum_{i=1}^N p_i \delta_{\mu_i}(x) dx = \sum_{i=1}^N p_i V(\mu_i).
\end{equation}
One can obtain the other two energies on the probability simplex via the same procedure. We list without derivation all the energy formula and corresponding gradient flows below:
\bequ\label{equ:all-func}
	\begin{aligned}
		\mcU\lp \rho \rp &= \int U\lp \rho\lp x \rp \rp dx = \sum_{i = 1}^N U\lp p_i \rp,		\\
   \dot{p_i} &= - \sqrt{p_i p_{i - 1} } \lp U'\lp p_i \rp - U'\lp p_{i - 1} \rp \rp + \sqrt{p_i p_{i + 1} } \lp U'\lp p_{i + 1} \rp - U'\lp p_{i} \rp \rp,\\
		\mcV\lp \rho \rp &= \int V\lp x \rp \rho(x)dx = \sum_{i = 1}^N V(\mu_i) p_i, 	\\
   \dot{p_i} & = - \sqrt{p_i p_{i - 1} } \lp V_i - V_{i - 1} \rp + \sqrt{p_i p_{i + 1} } \lp V_{i + 1} - V_{i} \rp,	\\
	 \mcW\lp \rho \rp &= \half \int\int W\lp x - y \rp \rho\lp x \rp\rho \lp y \rp dxdy = \half \sum_{i, j = 1}^N W(\mu_i - \mu_j) p_i p_j.		\\
\dot{p_i} & = - \sqrt{p_i p_{i - 1} } \lp \sum_{k = 1}^N W_{ik}p_k - \sum_{k = 1}^N W_{i-1, k}p_k \rp + \sqrt{p_i p_{i + 1} } \lp \sum_{k = 1}^N W_{i + 1, k}p_k - \sum_{k = 1}^N W_{ik}p_k \rp.
	\end{aligned}
\eequ
\subsubsection{Gradient flow of entropy in GMM}
In this subsection, we use the entropy functional which belongs to the internal functional as an examples and leave other derivations to interested readers. As illustrated in \Cref{general-gf}, the gradient flow of the entropy functional is nothing but the heat equation and it will also be the numerical scenario we investigated in \Cref{experiment}.

Recall the discrete entropy functional on a probability simplex is given by
\bequ
    H(p) = \sum_{i=1}^N p_i \log p_i,
\eequ
i.e. $U(p_i) = p_i \log p_i$. Differentiating to obtain $U'(p_i) = \log p_i + 1$ and using \cref{equ:all-func}, one obtains
\bequ\label{equ:1Dheat-gf}
    \dot{p}_i = -\sqrt{p_i p_{i - 1} } \log \frac{p_i}{p_{i-1}} + \sqrt{p_i p_{i + 1} } \log \frac{p_{i+1}}{p_{i}}.
\eequ

This parametric gradient flow equation can be interpreted from the perspective of the Markov transitional kernel. We use 2D dynamics to illustrate this point.
\begin{Ex}[2D dynamics]
For example, we write out the exact formula of entropy gradient flow for a probability simplex with three components
\bequn
	\lbb 
	\begin{aligned}
		\dot{p_1} = & \sqrt{p_1p_2}\log \frac{p_2}{p_1}	,		\\
		\dot{p_2} = & - \sqrt{p_1p_2}\log \frac{p_2}{p_1} + \sqrt{p_2p_3}\log \frac{p_3}{p_2},	\\
		\dot{p_3} = & - \sqrt{p_2p_3}\log \frac{p_3}{p_2}.
	\end{aligned}
	\rpt 
\eequn
We can also write out the evolution formula for this process in the matrix form, viewing the matrix as a Markovian kernel
\bequ
	\lp 
	\begin{aligned}
		& \dot{p_1}	\\
		& \dot{p_2}	\\
		& \dot{p_3}	\\
	\end{aligned}
	\rp = \begin{pmatrix}
		0  & \sqrt{\frac{p_1}{p_2}}\log \frac{p_2}{p_1} & 0		\\
		\sqrt{\frac{p_2}{p_1}}\log \frac{p_1}{p_2} & 0 & \sqrt{\frac{p_2}{p_3}}\log \frac{p_3}{p_2}	\\
		0 & \sqrt{\frac{p_3}{p_2}}\log \frac{p_2}{p_3} & 0
	\end{pmatrix}
	\lp 
	\begin{aligned}
		& p_1	\\
		& p_2	\\
		& p_3	\\
	\end{aligned}
	\rp .
\eequ
Notice that the Markovian kernel here depends on the density and thus is not temporally homogeneous. 

An observation shows that the unique equilibrium state for this Markov jump process is $p_1 = p_2 = p_3 = \frac{1}{3}$ and this fact actually holds for an arbitrary number of components, which is consistent with the continuous case.		\\
\indent
Introducing new parameter
\bequn
	a = \log \frac{p_1}{p_2}, b = \log \frac{p_2}{p_3},
\eequn
the Markovian kernel transform to
\bequn
	M_G = \begin{pmatrix}
		0  & \sqrt{\frac{p_1}{p_2}}\log \frac{p_2}{p_1} & 0		\\
		\sqrt{\frac{p_2}{p_1}}\log \frac{p_1}{p_2} & 0 & \sqrt{\frac{p_2}{p_3}}\log \frac{p_3}{p_2}	\\
		0 & \sqrt{\frac{p_3}{p_2}}\log \frac{p_2}{p_3} & 0
	\end{pmatrix} = \begin{pmatrix}
		0  & -ae^{\frac{a}{2}} & 0		\\
		ae^{-\frac{a}{2}} & 0 & -be^{\frac{b}{2}}	\\
		0 & be^{-\frac{b}{2}} & 0
	\end{pmatrix}.
\eequn
The evolution equation w.r.t. the parameter $a,b$ can be formulated as 
\bequn
	\lbb \begin{aligned}
		\dot{a} =&  \ \frac{d}{dt}{\log \frac{p_1}{p_2}} = \frac{\dot{p_1}}{p_1} - \frac{\dot{p_2}}{p_2} = -a \lp e^{\frac{a}{2}} + e^{- \frac{a}{2}} \rp - be^{-\frac{b}{2}},			\\
		\dot{b} =&  \frac{d}{dt}{\log \frac{p_2}{p_3}} = \frac{\dot{p_2}}{p_2} - \frac{\dot{p_3}}{p_3} = -b \lp e^{\frac{b}{2}} + e^{- \frac{b}{2}} \rp - ae^{\frac{a}{2}}.
	\end{aligned}
	\rpt
\eequn
\end{Ex}

Moreover, although our asymptotic analysis can not handle the situations when sample spaces are high dimensional, it is simple to generalize the gradient flow formula to higher dimensional cases. Here, we formally state the parametric gradient flow equation of 2-d heat flow under this framework without rigorous derivation. Consider a 2D grid $[0, N] \times [0, N]$ endowed with the natural graph structure which appears frequently as the underlying grid of finite difference method. Each node will have four neighbors, two of which in $x$-direction and others in $y$-direction. Then, we can simply write the RHS of 2D parametric heat equation as the sum of the RHS of two 1D parametric heat equation in x, y direction, i.e.
\bequ\label{equ:2Dheat-gf}
	\begin{aligned}
	\dot{p}_{ij} = & \ - \lp \sqrt{p_{i - 1,j}p_{i, ij}}\log \frac{p_{ij}}{p_{i - 1,j}} - \sqrt{p_{i + 1,j}p_{ij}}\log \frac{p_{i + 1,j}}{p_{ij}} \rpt			\\
	& \ \left. \qquad + \sqrt{p_{i,j - 1}p_{i, ij}}\log \frac{p_{ij}}{p_{i,j - 1}} - \sqrt{p_{i,j + 1}p_{ij}}\log \frac{p_{i,j + 1}}{p_{ij}} \rp.
	\end{aligned}
\eequ
We denote this as the 2D parametric heat equation and investigate it as a numerical scheme in \Cref{experiment}.

\subsubsection{Gradient flow of potential function in extended GMM}
In this subsection, we will derive step-by-step the gradient flow of a potential functional on the extended GMM \cref{equ:ex-W-metric}.

Recall that for an extended GMM there exists two classes of parameters, $\theta_i$s and $\mu_i$s. Therefore, we first take the gradient of the potential functional w.r.t. these parameters. In \cref{equ:all-func}, the discrete potential functional is written as
\bequ
    \mcV\lp \rho \rp = \sum_{i = 1}^N V(\mu_i) p_i.
\eequ
Notice the difference with homogeneous GMMs is that $V(\mu_i)$ is no longer a constant but rather depends on parameter $\mu_i$. Consequently, we have
\bequ
    \p_{p_i}\mcV = V(\mu_i), \quad \p_{{\mu}_i}\mcV = p_iV'(\mu_i).
\eequ
Now, recall that the WIM \cref{equ:ex-W-metric} is block-wise, one has
\bequ
    \begin{aligned}
        \dot{\theta_i} & \ = - \lp G_{\wtd W}^{ext} \rp_{\theta\theta, ii}^{-1} \p_{\theta_i} V\lp \theta, \mu \rp - \sum_j\lp G_{\wtd W}^{ext} \rp_{\theta\mu, ij}^{-1} \p_{\mu_j} V\lp \theta, \mu \rp,\\
        \dot{\mu_i} & \ = - \lp G_{\wtd W}^{ext} \rp_{\mu\mu, ii}^{-1} \p_{\mu_i} V\lp \theta, \mu \rp - \sum_j\lp G_{\wtd W}^{ext} \rp_{\mu\theta, ij}^{-1} \p_{\theta_j} V\lp \theta, \mu \rp,
    \end{aligned} 
\eequ
where the notation $\lp G_{\wtd W}^{ext} \rp_{\theta\theta, ii}^{-1}$ refers to the $(i, i)$ element of the block $\lp G_{\wtd W}^{ext} \rp_{\theta\theta}^{-1}$. Combine all the ingredients together, we conclude the following formula as the parametric gradient flow of potential functional on extended GMM:
\bequ\label{equ:ext-potential}
		\begin{aligned}
			\dot{\theta}_i = & \ - \frac{\sqrt{p_ip_{i + 1}}}{K\lp \sigma \rp}\lp \lp V\lp \mu_{i + 1} \rp - V\lp \mu_{i} \rp \rp - \frac{\mu_{i + 1} - \mu_i}{2}\lp  V'\lp \mu_i \rp + V'\lp \mu_{i + 1} \rp \rp \rp, \\
			& \ \quad i = 1,2, \cdots, N - 1,	\\
			\dot{\mu}_i = & \ - V'\lp \mu_i \rp + \frac{1}{K\lp \sigma \rp}\lp \sqrt{\frac{p_{i + 1}}{p_{i}}}\frac{\mu_{i + 1} - \mu_i}{2}\lp V\lp \mu_{i + 1} \rp - V\lp \mu_{i} \rp\rp \rpt	\\
			& \ \left. + \sqrt{\frac{p_{i - 1}}{p_{i}}}\frac{\mu_{i} - \mu_{i - 1}}{2}\lp V\lp \mu_{i} \rp - V\lp \mu_{i - 1} \rp\rp \rp, \quad i = 2, 3,\cdots, N - 1.
		\end{aligned}
\eequ
We will conduct some illustrative experiments on this model in the next section.

\section{Numerical experiments}\label{experiment}
Starting from parametric gradient flow equations derived in \cref{equ:gradient-flow-eq}, we will show that these equations can be viewed as numerical schemes for certain partial differential equation. Moreover, these schemes form sharp contrast with traditional method such as finite difference method, finite element method, and finite volume method in the sense that it origins from approximation using GMMs. We verifies the accuracy of these schemes using several partial differential equations.
\subsection{1D heat equation}
\begin{figure}[ht]
\label{1Dheat}
  \centering
  \centerline{\includegraphics[width=.8\linewidth]{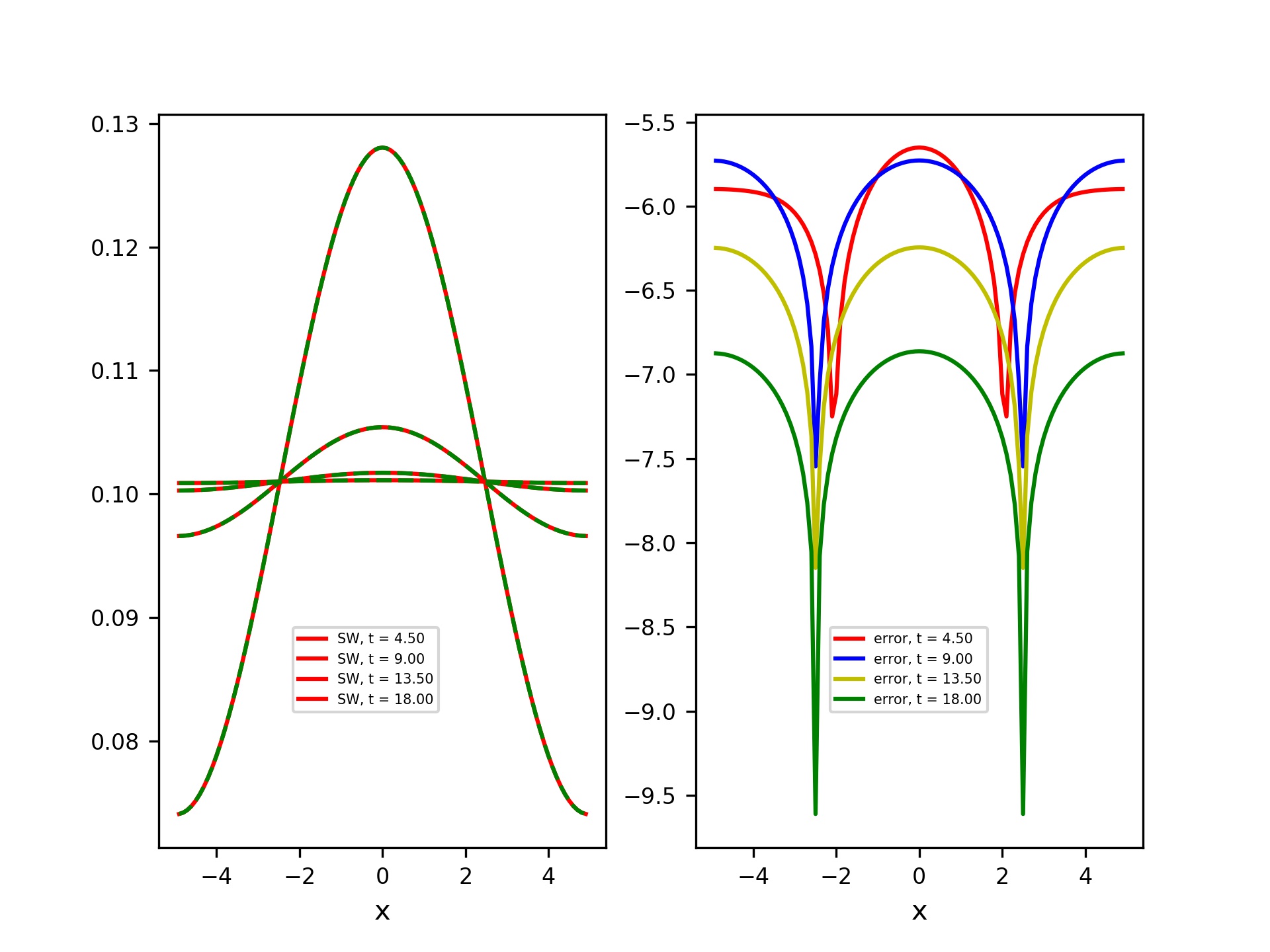}}
  \caption{This figure plots a simulation of the 1D heat flow via the discretization introduced in this paper. The parameters of the simulation are given by $\Delta x = 0.1, \Delta t = 0.001$. The initial distribution is given by $\rho\lp x; 0 \rp = \frac{e^{-\frac{x^2}{2}}}{\sqrt{2\pi}}$ and we use periodic boundary condition. In the left figure, we plot the solution using the scheme derived from scaling Wasserstein distance in red line and benchmark solution in green dashed line. The benchmark finite difference solution is solved with Crank-Nicolson scheme with the same parameters. In the right figure, we plot the log of the difference between these two solutions. It can be observed that the solution derived from the scaling Wasserstein metric has high accuracy.}
\end{figure}

We start from the 1D heat equation. Consider a GMM with all its means form a computational grid with gap $\Delta x$, i.e.
\bequ
    \mu_i = i\Delta x, \quad i = -n, -n+1, \cdots, n-1, n.
\eequ
Here we shift the indices to make the computational domain $[-5, 5] = [-n\Delta x, n\Delta x]$ symmetric w.r.t. $0$. We use periodic boundary condition on this domain and set the initial condition to be a Gaussian profile, i.e.
\bequ
    p_i = \frac{e^{-\frac{\mu_i^2}{2}}}{\sqrt{2\pi}}.
\eequ
Recall in \cref{equ:1Dheat-gf} we derive the parametric gradient flow equation of entropy functional, dividing RHS by spatial discretization size $\Delta x$ we obtain
\bequ\label{equ:heat-1Dscheme}
	\dot{p}_{i} = - \frac{1}{(\Delta x)^2} \lp \sqrt{p_{i - 1}p_i}\log \frac{p_i}{p_{i - 1}} - \sqrt{p_{i + 1}p_i}\log \frac{p_{i + 1}}{p_{i}} \rp.
\eequ
Now, if we think the $p_i$ as the density value on given grid point $\mu_i=i\Delta x$, i.e. $p_i = p(\mu_i)$, then the LHS of \cref{equ:heat-1Dscheme} can be thought of as a spatial discretization of the Laplacian operator. We implement this scheme and test it with Crank-Nicolson scheme to show its accuracy in \Cref{1Dheat}. From now on, we use SW (scaling Wasserstein) to represent the solution obtained via scaling Wasserstein metric and FD to denote the classical finite difference solution.

\subsection{2D heat equation}
\begin{figure}[ht]
\label{2Dheat}
  \centering
\centerline{\includegraphics[width=.8\linewidth]{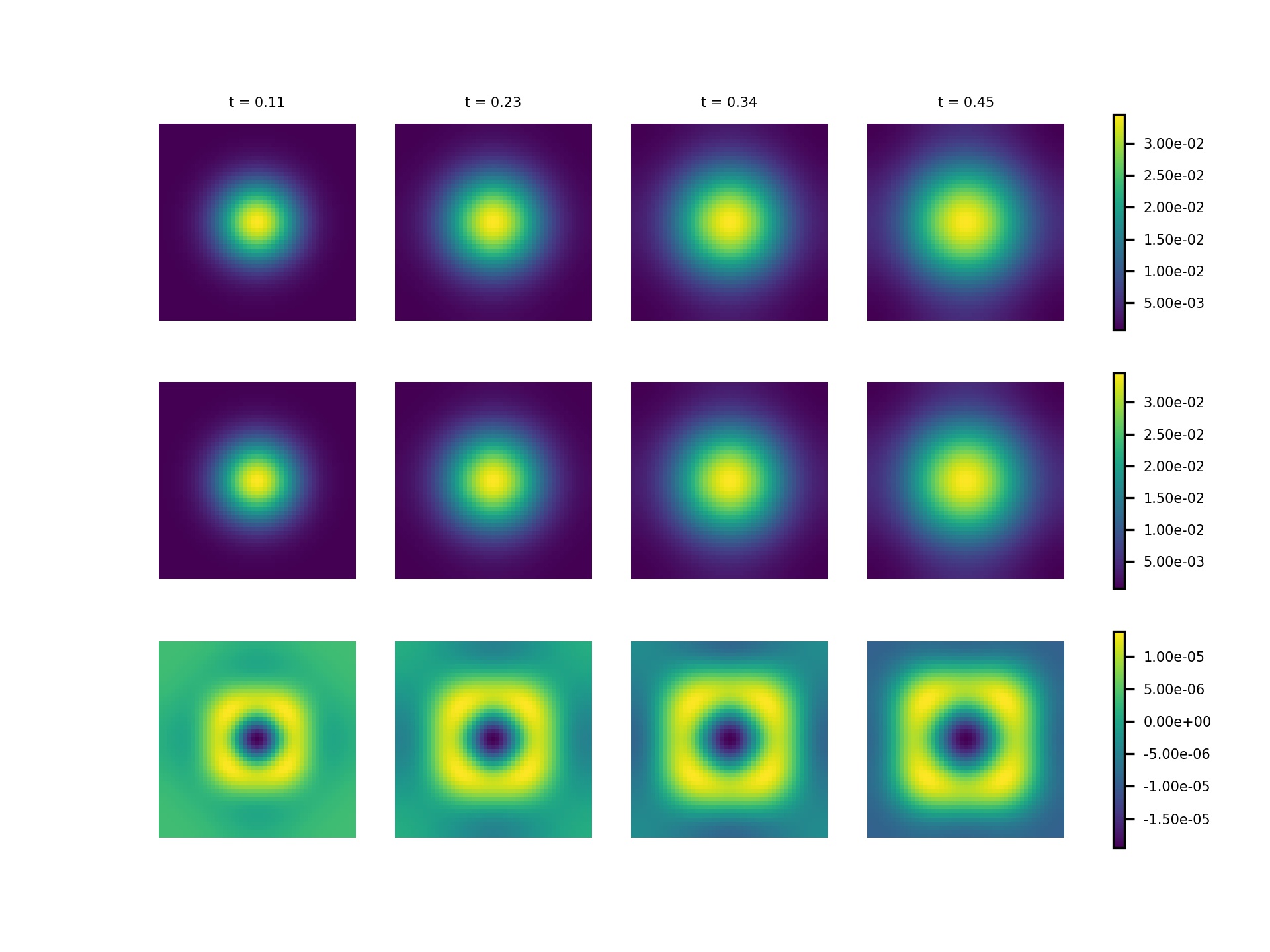}}
  \caption{This figure plots a simulation of the 2D heat flow via the discretization introduced in this paper. The parameters of the simulation are given by $\Delta x = 0.1, \Delta t = 0.001$. The initial distribution is given by $\rho\lp x; 0 \rp = \frac{e^{-\frac{x^2+y^2}{2}}}{2\pi}$ and we use periodic boundary condition. In (a), we plot the snapshot of the density obtained via SW at different time step. In (b), we plot the density obtained via FD at the same time steps. In (c), we plot the difference between these two solutions. It can be observed that the solution derived from scaling Wasserstein metric has high accuracy.}
\end{figure}

Similarly, recall in \cref{equ:2Dheat-gf} we formally derive the parametric gradient flow equation of entropy functional in 2D. Again we divide the RHS by the spatial discretization size $\Delta x$ to obtain
\bequn\label{heat-2Dscheme}
	\begin{aligned}
	\dot{p}_{ij} = & \ - \frac{1}{(\Delta x)^2} \lp \sqrt{p_{i - 1,j}p_{i, ij}}\log \frac{p_{ij}}{p_{i - 1,j}} - \sqrt{p_{i + 1,j}p_{ij}}\log \frac{p_{i + 1,j}}{p_{ij}} \rpt			\\
	& \ \left. \qquad\qquad + \sqrt{p_{i,j - 1}p_{i, ij}}\log \frac{p_{ij}}{p_{i,j - 1}} - \sqrt{p_{i,j + 1}p_{ij}}\log \frac{p_{i,j + 1}}{p_{ij}} \rp.
	\end{aligned}
\eequn
Comparing with \cref{equ:heat-1Dscheme}, the 2D scheme can be viewed as a summation of two 1D schemes on $x, y$ directions, which manifests the fact that the Laplacian is separable, i.e. the 2D Laplacian is simply the summation of two 1D Laplaces. Simulating \Cref{heat-2Dscheme} on a periodic domain with Gaussian initial distribution is presented in \Cref{2Dheat}. The parameters of the simulation are given by $\Delta x = 0.1, \Delta t = 0.001$. The initial distribution is given by $\rho\lp x; 0 \rp = e^{-\frac{x^2+y^2}{2}}/2\pi$. In (a), we plot the snapshot of the density obtained via SW at different time step. In (b), we plot the density obtained via FD at the same time steps. In (c), we plot the difference between these two solutions. It can be observed that the solution derived from scaling Wasserstein metric has high accuracy.
\subsection{1D transport equation on extended GMM}
\begin{figure}[ht]
  \centering
    \label{eGMM-potential}
\centerline{\includegraphics[width=\linewidth]{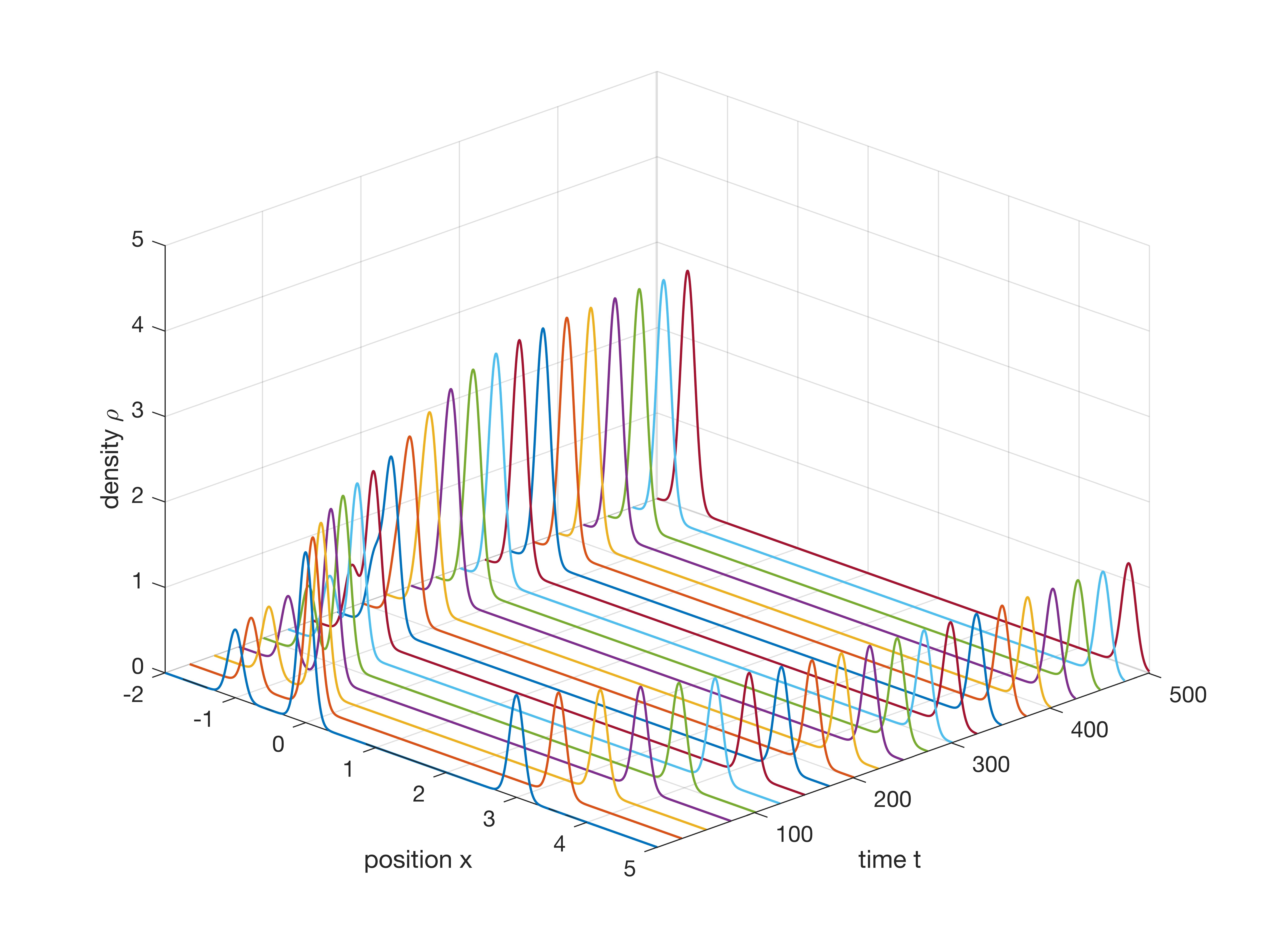}}
  \caption{In this figure, we plot the simulation of the Wasserstein gradient flow of the potential function $V\lp x \rp = \sin x$ in the extended GMM \Cref{sec:ext-GMM}. The initial data is given by $\rho_0 = 0.2 * \mcN\lp -1, 0.1 \rp + 0.5 * \mcN\lp 0, 0.1 \rp + 0.3 * \mcN\lp 3, 0.1 \rp$. The time step is set as $\Delta t = 0.01$. The iteration epoch is set as $5000$. Only first two components collapse to a common minimizer $x = - \frac{\pi}{2}$ while the rightmost component converges to $x = \frac{3\pi}{2}$. We observe the mode degeneracy of two means $\mu_1, \mu_2$ in the above figure.}
\end{figure}

In this last numerical example, we present some preliminary experiments on the extended GMM \Cref{sec:ext-GMM}. We simulate a Wasserstein gradient flow of the potential function $V\lp x \rp = \sin x$, i.e.
\bequ
    \p_t \rho(x) - \nabla \cdot (\rho(x) \nabla V(x)) = 0.
\eequ
The discretization according to the extended GMM is presented in \cref{equ:ext-potential}. We use forward Euler for the temporal discretization. The initial data is given by $\rho_0 = 0.2 * \mcN\lp -1, 0.1 \rp + 0.5 * \mcN\lp 0, 0.1 \rp + 0.3 * \mcN\lp 3, 0.1 \rp$. The time step is set as $\Delta t = 0.01$ with number of iteration $5000$. As the minimums of $V(x)$ appear in period $2\pi$, each Gaussian component will be driven to the nearest minimum by the gradient flow. Therefore, the first two components collapse to a common minimum $x = - \frac{\pi}{2}$ while the rightmost component converges to $x = \frac{3\pi}{2}$. We observe the mode degeneracy of two means $\mu_1, \mu_2$ in \Cref{eGMM-potential}.

\section{Discussion}\label{discussion}
In this paper, we study the Wasserstein pullback metric on the Gaussian mixture models. Our key observation is that as the variance of Gaussian components tends to $0$, the scaling WIM reveals a diagonal structure. Based on this, we construct and prove the convergence of a metric structure under the scaling limit that variance tends to $0$. We establish its analytic form via asymptotic analysis and introduce other related models with inhomogeneous gaps, higher-order terms, and extra degrees of freedom for the mean value parameters. Next, the gradient flows on this metric space are introduced, focusing on their relations with the numerical schemes of partial differential equations. We verify the correctness and accuracy of these schemes on 1D and 2D heat equations and present some preliminary experiments on the extended Gaussian mixture models.

A systematic study of the scaling Wasserstein geometry is needed. Questions such as functional inequalities and geometric properties remain to be explored. Connections to spectral graph theory deserve to be studied. Besides, one also needs to establish a correspondent theory for high-dimensional sample spaces. The current method to obtain the conclusion on 1D sample spaces relies significantly on the closed-form solution of WIMs, which is often absent in high dimensions. Lastly, apart from gradient flows, Hamiltonian flows and related mathematical physics equations, including Schr{\"o}dinger equations, Schr{\"o}dinger bridge problems, etc., can also be studied in the Gaussian mixture models. Detailed convergence analysis of WIM from GMM models to the Wasserstein metric in density space is left in future work. 

\nocite{khan2022optimal}

\bibliographystyle{abbrv}
\bibliography{Scaling}

\begin{thebibliography}{10}

\bibitem{NG}
S.~Amari.
\newblock Natural {{Gradient Works Efficiently}} in {{Learning}}.
\newblock {\em Neural Computation}, 10(2):251--276, 1998.

\bibitem{IG}
S.~Amari.
\newblock {\em Information Geometry and Its Applications}.
\newblock Number volume 194 in Applied mathematical sciences. {Springer},
  Japan, 2016.

\bibitem{Ambrosio2008Gradient}
L.~Ambrosio, N.~Gigli, and {Savar\'e Giuseppe}.
\newblock {\em Gradient {{Flows}}: In {{Metric Spaces}} and in the {{Space}} of
  {{Probability Measures}}}.
\newblock {Birkh\"auser Basel}, Basel, 2005.

\bibitem{ay2017information}
N.~Ay, J.~Jost, H.~L{\^e}, and L.~Schwachh{\"o}fer.
\newblock {\em Information Geometry}.
\newblock Ergebnisse der Mathematik und ihrer Grenzgebiete. 3. Folge / A Series
  of Modern Surveys in Mathematics. Springer International Publishing, 2017.

\bibitem{chen2018optimal}
Y.~Chen, T.~T. Georgiou, and A.~Tannenbaum.
\newblock Optimal transport for gaussian mixture models.
\newblock {\em IEEE Access}, 7:6269--6278, 2018.

\bibitem{NIPS2013_af21d0c9}
M.~Cuturi.
\newblock Sinkhorn distances: Lightspeed computation of optimal transport.
\newblock In C.~Burges, L.~Bottou, M.~Welling, Z.~Ghahramani, and
  K.~Weinberger, editors, {\em Advances in Neural Information Processing
  Systems}, volume~26. Curran Associates, Inc., 2013.

\bibitem{e25050786}
Q.~Feng and W.~Li.
\newblock Entropy dissipation for degenerate stochastic differential equations
  via sub-riemannian density manifold.
\newblock {\em Entropy}, 25(5), 2023.

\bibitem{refId0}
{Gangbo, Wilfrid}, {Li, Wuchen}, and {Mou, Chenchen}.
\newblock Geodesics of minimal length in the set of probability measures on
  graphs.
\newblock {\em ESAIM: COCV}, 25:78, 2019.

\bibitem{gao2022master}
Y.~Gao, W.~Li, and J.-G. Liu.
\newblock Master equations for finite state mean field games with nonlinear
  activations, 2022.

\bibitem{doi:10.1137/19M1243440}
P.~Gladbach, E.~Kopfer, and J.~Maas.
\newblock Scaling limits of discrete optimal transport.
\newblock {\em SIAM Journal on Mathematical Analysis}, 52(3):2759--2802, 2020.

\bibitem{GLADBACH2020204}
P.~Gladbach, E.~Kopfer, J.~Maas, and L.~Portinale.
\newblock Homogenisation of one-dimensional discrete optimal transport.
\newblock {\em Journal de Mathématiques Pures et Appliquées}, 139:204--234,
  2020.

\bibitem{191765}
D.~(https://math.stackexchange.com/users/31254/donantonio).
\newblock Help with integrating $\int_0^{\infty} \frac{(\log x)^2}{x^2 + 1}
  \operatorname d\!x$ - contour integration?
\newblock Mathematics Stack Exchange.
\newblock URL:https://math.stackexchange.com/q/191765 (version: 2017-05-02).

\bibitem{RePEc:spr:aistmt:v:74:y:2022:i:1:d:10.1007_s10463-021-00788-1}
S.~ichi Amari and T.~Matsuda.
\newblock {Wasserstein statistics in one-dimensional location scale models}.
\newblock {\em Annals of the Institute of Statistical Mathematics},
  74(1):33--47, February 2022.

\bibitem{JKO}
R.~Jordan, D.~Kinderlehrer, and F.~Otto.
\newblock The {{Variational Formulation}} of the {{Fokker}}--{{Planck
  Equation}}.
\newblock {\em SIAM Journal on Mathematical Analysis}, 29(1):1--17, 1998.

\bibitem{khan2022optimal}
G.~Khan and J.~Zhang.
\newblock When optimal transport meets information geometry.
\newblock {\em Information Geometry}, 5(1):47--78, 2022.

\bibitem{kriegl1997convenient}
A.~Kriegl and P.~W. Michor.
\newblock {\em The convenient setting of global analysis}, volume~53.
\newblock American Mathematical Soc., 1997.

\bibitem{LiG}
W.~Li.
\newblock Geometry of probability simplex via optimal transport.
\newblock {\em arXiv:1803.06360 [math]}, 2018.

\bibitem{10.1063/5.0012605}
W.~Li.
\newblock {Hessian metric via transport information geometry}.
\newblock {\em Journal of Mathematical Physics}, 62(3):033301, 03 2021.

\bibitem{bregmanTIG}
W.~Li.
\newblock Transport information bregman divergences.
\newblock {\em Information Geometry}, 4(2):435--470, 2021.

\bibitem{li2022mean}
W.~Li and L.~Lu.
\newblock Mean field information hessian matrices on graphs, 2022.

\bibitem{LI2022109645}
W.~Li and F.~Rubio.
\newblock On a prior based on the wasserstein information matrix.
\newblock {\em Statistics \& Probability Letters}, 190:109645, 2022.

\bibitem{li2023wasserstein}
W.~Li and J.~Zhao.
\newblock Wasserstein information matrix.
\newblock {\em Information Geometry}, pages 1--53, 2023.

\bibitem{Lottgeo}
J.~Lott.
\newblock Some geometric calculations on wasserstein space.
\newblock {\em Communications in Mathematical Physics}, 277(2):423--437, 2008.

\bibitem{LV}
J.~Lott and C.~Villani.
\newblock Ricci curvature for metric-measure spaces via optimal transport.
\newblock {\em Annals of Mathematics}, Vol. 169(No. 3):903--991, 2009.

\bibitem{MAAS20112250}
J.~Maas.
\newblock Gradient flows of the entropy for finite markov chains.
\newblock {\em Journal of Functional Analysis}, 261(8):2250--2292, 2011.

\bibitem{Malago2}
L.~Malag\`o and G.~Pistone.
\newblock Natural {{Gradient Flow}} in the {{Mixture Geometry}} of a {{Discrete
  Exponential Family}}.
\newblock {\em Entropy}, 17(12):4215--4254, 2015.

\bibitem{Mielke_2011}
A.~Mielke.
\newblock A gradient structure for reaction–diffusion systems and for
  energy-drift-diffusion systems*.
\newblock {\em Nonlinearity}, 24(4):1329, mar 2011.

\bibitem{murray2012asymptotic}
J.~D. Murray.
\newblock {\em Asymptotic analysis}, volume~48.
\newblock Springer Science \& Business Media, 2012.

\bibitem{otto2000generalization}
F.~Otto and C.~Villani.
\newblock Generalization of an inequality by talagrand and links with the
  logarithmic sobolev inequality.
\newblock {\em Journal of Functional Analysis}, 173(2):361--400, 2000.

\bibitem{peyre2019computational}
G.~Peyr{\'e}, M.~Cuturi, et~al.
\newblock Computational optimal transport: With applications to data science.
\newblock {\em Foundations and Trends{\textregistered} in Machine Learning},
  11(5-6):355--607, 2019.

\bibitem{reynolds2009gaussian}
D.~A. Reynolds et~al.
\newblock Gaussian mixture models.
\newblock {\em Encyclopedia of biometrics}, 741(659-663), 2009.

\bibitem{slepcev2022nonlocal}
D.~Slepcev and A.~Warren.
\newblock Nonlocal wasserstein distance: Metric and asymptotic properties,
  2022.

\bibitem{villani2021topics}
C.~Villani.
\newblock {\em Topics in optimal transportation}, volume~58.
\newblock American Mathematical Soc., 2021.

\bibitem{doi:10.1002/cpa.20060}
M.-K. von Renesse and K.-T. Sturm.
\newblock Transport inequalities, gradient estimates, entropy and ricci
  curvature.
\newblock {\em Communications on Pure and Applied Mathematics}, 58(7):923--940,
  2005.

\bibitem{Wangflow}
Y.~Wang and W.~Li.
\newblock Accelerated information gradient flow.
\newblock {\em Journal of Scientific Computing}, 90(1):11, 2021.

\end{thebibliography}

\end{document}